\UseRawInputEncoding
\documentclass[a4wide,english]{amsart}

\usepackage{amsmath}

\usepackage[all]{xy}
\CompileMatrices
\newdir{ >}{{}*!/-10pt/\dir{>}}

\theoremstyle{plain}
\newtheorem{theorem}{Theorem}[section]

\theoremstyle{definition}

\newtheorem{emptythm}[theorem]{}

\newcommand{\OX}{{\mathcal O}_{X}}

\newcommand{\OXx}{{\mathcal O}_{X,x}}

\newcommand{\HOM}{{\mathcal H}om}

\newcommand{\Az}{\sheaf{A}}

\newcommand{\dual}[1]{{\mathfrak D}^{#1}}

\newcommand{\cod}[1]{\mu_{#1}}

\newcommand{\inv}{\tau}
\newcommand{\binv}{\nu}

\newcommand{\cw}[2]{\tilde{\W}^{#2}(#1)}
\newcommand{\cws}[3]{\tilde{W}^{#3}_{#2}(#1)}

\newcommand{\too}{\longrightarrow}
\newcommand{\sheaf}[1]{{\mathcal #1}}
\newcommand{\smb}{{\scriptscriptstyle \bullet}}

\newcommand{\ol}[1]{\overline{#1}}

\DeclareMathOperator{\hgt}{ht}
\DeclareMathOperator{\Hom}{Hom}

\DeclareMathOperator{\Ext}{Ext}

\DeclareMathOperator{\Ker}{Ker}

\DeclareMathOperator{\id}{id}

\DeclareMathOperator{\Spec}{Spec}

\DeclareMathOperator{\supp}{supp}

\DeclareMathOperator{\codimension}{codim}

\DeclareMathOperator{\Zar}{\mathrm{Zar}}

\DeclareMathOperator{\ceZ}{\mathrm{Z}}

\DeclareMathOperator{\W}{\mathrm{W}}
\DeclareMathOperator{\GW}{\mathrm{GW}}

\DeclareMathOperator{\tr}{\mathrm{tr}}

\newcommand{\Mcat}[1]{\sheaf{M}(#1)}
\newcommand{\Mqc}[1]{\sheaf{M}_{qc}(#1)}
\newcommand{\Pcat}[1]{\sheaf{P}(#1)}

\newcommand{\Mfg}[1]{\sheaf{M}_{c}(#1)}
\newcommand{\MfgS}[2]{\sheaf{M}_{c\, ,#2}(#1)}
\newcommand{\MqcS}[2]{\sheaf{M}_{qc\, ,#2}(#1)}

\newcommand{\DP}[1]{\mathrm{D}^{b}(\Pcat{#1})}
\newcommand{\DPS}[2]{\mathrm{D}_{#2}^{b}(\Pcat{#1})}
\newcommand{\DPF}[2]{\mathrm{D}^{b}(\Pcat{#1})^{(#2)}}

\newcommand{\Dc}[1]{\mathrm{D}^{b}_{c}(\Mqc{#1})}

\newcommand{\DMqc}[1]{\mathrm{D}^{b}(\Mqc{#1})}

\newcommand{\DcS}[2]{\mathrm{D}^{b}_{c\, ,#2}(\Mqc{#1})}
\newcommand{\DcSF}[3]{\mathrm{D}^{b}_{c\, ,#2}(\Mqc{#1})^{(#3)}}

\newcommand{\Dfg}[1]{\mathrm{D}^{b}_{c}(\Mqc{#1})}

\newcommand{\DcF}[3]{\mathrm{D}^{b}_{c}(\Mcat{#1})^{(#2)}_{#3}}
\newcommand{\DfgF}[2]{\mathrm{D}^{b}_{c}(\Mqc{#1})^{(#2)}}

\DeclareMathOperator{\HM}{\mathrm{H}}

\newcommand{\Z}{\mathbb{Z}}

\newcommand{\ideal}{\mathfrak{a}}
\newcommand{\Jideal}{\mathfrak{J}}
\newcommand{\maxid}{\mathfrak{m}}

\newcommand{\cf}{\textsl{cf.}\ }
\newcommand{\eg}{\textsl{e.g.}\ }
\newcommand{\ie}{\textsl{i.e.}\ }

\begin{document}

\title[On the Gersten conjecture for hermitian Witt groups]
{On the Gersten conjecture for hermitian Witt groups}

\author{Stefan Gille}
\email{gille@ualberta.ca}
\address{Department of Mathematical and Statistical Sciences,
University of Alberta, Edmonton T6G 2G1, Canada}

\author{Ivan Panin}
\email{paniniv@pdmi.ras.ru}
\address{Petersburg Department of Steklov Institute of Mathematics,
27, Fontanka 191011, St.Petersburg, Russia}

\thanks{The work of S.\ G.\ has been supported by an NSERC grant.}


\date{January 24, 2022}

\begin{abstract}
We prove that the hermitian Gersten-Witt complex is exact for Azumaya algebras with involution
of the first- or second kind over a regular local ring, which is essentially smooth over a field, or
over a discrete valuation ring.
\end{abstract}

\maketitle

\bigbreak

\section{Introdcution}
\label{IntroSect}\bigbreak

\noindent
Let~$R$ be a regular integral domain of finite Krull dimension with fraction field~$K$ of characteristic not two,
and $(A,\inv)$ an Azumaya algebra with involution of the first- or second kind over~$R$. In~\cite{Gi07,Gi09,Gi13}
the first named author has constructed a complex, the so called $\epsilon$-hermitian Gersten-Witt complex of~$(A,\inv)$,
$\epsilon\in\{\pm 1\}$:

{\footnotesize
$$
0\too\W_{\epsilon}(A,\inv)\too\W_{\epsilon}(K\otimes_{R}(A,\inv))\too
\bigoplus\limits_{\hgt q=1}\W_{\epsilon}(k(q)\otimes_{R}(A,\inv))\too\ldots\qquad\qquad\qquad
$$
$$
\qquad\qquad\qquad\qquad\qquad\qquad\qquad\qquad\ldots\too
\bigoplus\limits_{\hgt q=\dim R}\W_{\epsilon}(k(q)\otimes_{R}(A,\inv))\too 0\, ,
$$
}

\noindent
where $k(q)$ denotes the residue field of $q\in\Spec R$, and
$\W_{\epsilon}(k(q)\otimes_{R}(A,\inv))$, $\epsilon\in\{\pm 1\}$, denotes the $\epsilon$-hermitian
Witt group of the central simple $k(q)$-algebra $k(q)\otimes_{R}A$ with involution $\id_{k(q)}\otimes\,\inv$.
This construction is the natural generalization of the one of Balmer and Walter~\cite{BaWa02}
for Witt groups of symmetric forms.

\smallbreak

The Gersten conjecture claims that if~$R$ is a regular local ring then this complex is exact. In the symmetric
case, \ie $(A,\inv)=(R,\id_{R})$, this conjecture has been verified in many instances, \eg for regular local rings
of dimension~$\leq 4$ by Balmer and Walter~\cite{BaWa02}, or for regular local rings~$R$ which contain a field
(of characteristic not~$2$), by Balmer, Walter, and the authors~\cite{BaGiPaWa02}. In the hermitian case the
first named author has given a proof if~$R$ is regular local and essentially smooth over a field, and~$(A,\inv)$ 
is extended from the base field in~\cite{Gi09,Gi13}. These papers claim the conjecture also in the non constant case,
\ie if~$(A,\inv)$ is not coming from the base field, but the proof is flawed, see our Remark~~\ref{wrongGWPfRem}
(we fix this gap here). Recently Bayer-Fluckiger, First, and Parimala~\cite{BaFiPa19} have verified the conjecture if
$\dim R\leq 2$, and if~$\dim R\leq 4$ and~$A$ is of odd index.

\smallbreak

In this article we prove the conjecture for regular local rings, which are essential smooth over a discrete valuation ring.

\smallbreak

\noindent
{\bf Theorem.}
{\it
Let~$R$ be a integral domain, which is smooth over a discrete valuation ring, or over a field, $\tilde{R}$ a localization
of~$R$ at a prime ideal, and $(\tilde{A},\tilde{\inv})$ an Azumaya algebra with involution of the first-
or second kind over~$\tilde{R}$. Then the Gersten conjecture holds for~$(\tilde{A},\tilde{\inv})$.
}

\smallbreak

\noindent
Note that this theorem is new even in the symmetric case, \ie $(\tilde{A},\tilde{\inv})=(\tilde{R},\id_{\tilde{R}})$.
Using Popescu's desingularization theorem~\cite{Po85,Po86} our result implies the conjecture also for
Azumaya algebras with involution over a regular local ring, which either contains a field, or which is geometrically
regular over a discrete valuation ring.

\smallbreak

We give now a short sketch of the proof of the main theorem, which
is in its essence an adaption of Quillen's~\cite{Qu73} proof of the Gersten conjecture
in $K$-theory to hermitian Witt groups.

\smallbreak

By assumption we have $\tilde{R}=R_{P}$ for some
prime ideal~$P$ of~$R$, and replacing~$R$ by a localization we can assume that
$(\tilde{A},\tilde{\inv})=\tilde{R}\otimes_{R}(A,\inv)$ for some Azumaya algebra
with involution~$(A,\inv)$ over~$R$. Denote by~$\Dc{\tilde{A}}$ the bounded derived
category of complexes of $\tilde{A}$-modules with finitely generated homology modules,
and by~$\DfgF{\tilde{A}}{p}$, $p\geq 0$~an integer, the full subcategory consisting of
complexes~$M_{\smb}$ with $\codimension_{\Spec R}\supp M_{\smb}\geq p$. A finite
injective resolution of~$\tilde{R}$ considered as an element in the bounded derived
category~$\Dc{\tilde{R}}$ is a dualizing complex and so induces a duality on~$\Dc{\tilde{A}}$
as well as on $\DfgF{\tilde{A}}{p}$ giving these categories the structure of triangulated categories
with duality in the sense of Balmer~\cite{Ba00}.

\smallbreak

By construction the Gersten conjecture is equivalent to the assertion that
the natural functor $\DfgF{A}{p+1}\too\DfgF{A}{p}$ induces the zero map on the associated
triangular Witt groups for all $p\geq 0$. In Section~\ref{GWConjSmoothDVRSect} we show
that this follows in turn from the following result (see Lemma~\ref{GWPfmaintechLem}
for a precise formulation including in particular the involved dualizing complexes):

\smallbreak

\noindent
{\it
Let~$t\in R$ be a non zero divisor, such that $R':=R/Rt$ is flat over the base ring,
$\pi:R\too R'$ the quotient map, and $\gamma:R\too\tilde{R}=R_{P}$ the localization
morphism. Then
$$
\gamma^{\ast}\tr_{\pi}(x)\, =\, 0\quad\mbox{in~$W^{i}(\DfgF{\tilde{A}}{p})$}
$$
for all $x\in\W^{i}(\DfgF{R'\otimes_{R}A}{p})$.
}

\smallbreak

\noindent
Here~$\tr_{\pi}$ stands for the transfer map along~$\pi$ and~$\W^{i}$ for
the $i$th triangular Witt group.

\smallbreak

As this is merely an outline of the idea of proof we do not mention for simplicity here and
in the following the involved dualities, see Section~\ref{PfGWPfmainLemSect} for this.

\smallbreak

The main geometric ingredient in the proof of above claim is the normalization lemma of Quillen~\cite{Qu73},
respectively its generalization by Gillet and Levine~\cite{GiLe87} in case the base ring is a discrete
valuation ring. This result coupled with Zariski's main theorem provides us with a commutative
diagram
$$
\xymatrix{
 & & & R'_{P}
\\
 & C' \ar[urr]^-{s} & & R' \ar[u]_-{\gamma'} \ar[ll]^-{u}
\\
D \ar[ur]^-{\alpha'} & & &
\\
 & \tilde{R} \ar[ul]^-{\delta} & & R \ar[uu]_-{\pi} \ar[ll]^-{\gamma} \rlap{\, ,} 
}
$$
where~$u$ is essentially smooth (and so~$C'$ is Gorenstein), $\gamma'$ is the localization
morphism, $s$~a regular immersion of codimension one, $\delta$~finite, the by~$\alpha'$ induced
morphism $\Spec C'\too\Spec D$ an open immersion, and $s\circ\alpha'$ is surjective.

\smallbreak

Set $(A',\inv'):=R'\otimes (A,\inv)$. In general the Azumaya algebras with involution $u^{\ast}(A',\inv')$
and $(\alpha'\circ\delta)^{\ast}(\tilde{A},\tilde{\inv})$ are not isomorphic. In particular, there are two
transfer maps along $s:C'\too R'_{P}$:
$$
\tr^{1}_{s}\, :\;\W^{i}(\DfgF{\tilde{A}/\tilde{A}t}{p})\,\too\,\W^{i}(\DfgF{u^{\ast}(A)}{p})
$$
and
$$
\tr^{2}_{s}\, :\;\W^{i}(\DfgF{\tilde{A}/\tilde{A}t}{p})\,\too\,\W^{i}(\DfgF{(\alpha'\circ\delta)^{\ast}(\tilde{A})}{p})\, .
$$
Now by the zero theorem for the transfer~\cite[Thm.\ 6.3]{Gi13} we have
$\tr^{1}_{s}(\gamma'^{\ast}(x))=0$, and using an excision lemma we show that
$$
\gamma^{\ast}(\tr_{\pi}(x))\, =\,\tr_{\delta}\Big[\, (\alpha'^{\ast})^{-1}\big(\tr^{2}_{s}(\gamma'^{\ast}(x))\big)\,\Big]\, .
$$
Hence if $u^{\ast}(A',\inv')\simeq (\alpha'\circ\delta)^{\ast}(\tilde{A},\tilde{\inv})$, which is for instance the
case if~$(\tilde{A},\tilde{\inv})$ is extended from the base ring, this concludes the proof. However -- as
already mentioned -- these algebras with involutions are in general not isomorphic. In this case we remedy
this obstruction using a construction of Ojanguren and the second named author~\cite{OjPa01}: There exists
a smooth morphism of relative dimension zero $\kappa:C'\too\tilde{C}$, such that
$$
\kappa^{\ast}\big(\, u^{\ast}(A',\inv')\,\big)\,\simeq\,
\kappa^{\ast}\big(\,(\alpha'\circ\delta)^{\ast}(\tilde{A},\tilde{\inv})\,\big)\, ,
$$
and satisfying another technical property, which is crucial since it implies that above
morphism $s:C'\too R'_{P}$ factors via~$\kappa$ and a regular immersion $\beta':\tilde{C}\too R'_{P}$.

\medbreak

This is done in the last Section~\ref{PfGWPfmainLemSect} of the paper.
The content of the rest of the article is as follows. In Sections~\ref{WittTheorySect} and~\ref{AzInvSect}
we recall the basic definitions of triangulated and derived (hermitian) Witt theory. Section~\ref{DualComplSect}
fixes some notations and recalls dualizing complexes. The following Section~\ref{CohWGrSect} is a jog
through coherent hermitian Witt theory of algebras with involutions over (commutative) rings with
dualizing complexes.

\smallbreak

In Section~\ref{techLemSect} we prove the above mentioned excision lemma for the
transfer (for simplicity only in the special situation we use it), and in Section~\ref{GWComplSect}
we recall the construction of the hermitian Gersten-Witt complex as a well as the formulation of the Gersten
conjecture and two of its consequences.

\bigbreak

\goodbreak
\section{Review of Witt theory of categories with duality}
\label{WittTheorySect}\bigbreak

\begin{emptythm}
\label{ExCatDSubSect}
{\bf Exact categories with duality.}
Throughout this work we assume that the $\Hom$-groups of additive
categories are uniquely $2$-divisible. In particular we assume that
schemes have $1/2$ in their global sections.

\smallbreak

An {\it exact category with duality} is a triple
$(\sheaf{E},\vee,\varpi)$, where~$\vee:\sheaf{E}\too\sheaf{E}$ is a contravariant exact
functor and~$\varpi$ a natural isomorphism $\id_{\sheaf{E}}\xrightarrow{\simeq}\vee\circ\vee$ satisfying
$\varpi_{M}^{\vee}=\varpi^{-1}_{M^{\vee}}$ for all $M\in\sheaf{E}$. 

\smallbreak

A {\it $\epsilon$-symmetric space}, $\epsilon\in\{\pm 1\}$, is a pair $(M,\varphi)$, where $\varphi:M\too M^{\vee}$
is an isomorphism in~$\sheaf{E}$, such that $\varphi^{\vee}\circ\varpi_{M}=\epsilon\cdot\varphi$.
Two $\epsilon$-symmetric spaces $(M_{1},\varphi_{1})$ and $(M_{2},\varphi_{2})$ are
called {\it isometric} if there exists an isomorphism $\theta:M_{1}\xrightarrow{\simeq}M_{2}$,
such that $\varphi_{1}=\theta^{\vee}\circ\varphi_{2}\circ\theta$. The associated {\it Witt group of
$\epsilon$-symmetric spaces} will be denoted $\W_{\epsilon}(\sheaf{E},\vee)$.
This is the Grothendieck group of the isometry classes of $\epsilon$-symmetric spaces
with the orthogonal sum as addition modulo the so called {\it metabolic spaces}.
\end{emptythm}

\begin{emptythm}
\label{TrCatDSubSect}
{\bf Triangulated categories with duality.}
We refer to the works~\cite{Ba00,Ba01a} of Balmer for details and more information.

\smallbreak

A {\it triangulated category with $\delta$-exact duality}, $\delta\in\{\pm 1\}$, is a triple $(\sheaf{T},\vee,\varpi)$
consisting of a triangulated category~$\sheaf{T}$, a $\delta$-exact duality $\vee$, and an isomorphism~$\varpi$
to the bidual satisfying the same axioms as the one for exact categories with duality. If the isomorphism to
the bidual~$\varpi$ is clear from the context we also say that the pair $(\sheaf{T},\vee)$, is a triangulated
category with duality.

\smallbreak

Denote by~$T$ the translation functor of~$\sheaf{T}$. Then $T^{i}\circ\vee$ is a
$(-1)^{i}\delta$-exact duality, and $(\sheaf{T},\, T^{i}\circ\vee,\, (-1)^{\frac{i(i+1)}{2}}\varpi)$
is a triangulated category with duality. A {\it $i$-symmetric space} is a symmetric space in
$(\sheaf{T},T^{i}\circ\vee,(-1)^{\frac{i(i+1)}{2}}\varpi)$. Isometry and the orthogonal sum of
spaces are defined as for exact categories with duality. The {\it $i$th triangular Witt group} of
$(\sheaf{T},\vee,\varpi)$, denoted $\W^{i}(\sheaf{T},\vee,\varpi)$, $i\in\Z$, or $\W^{i}(\sheaf{T},\vee)$,
respectively, is the Grothendieck-Witt group of the isometry classes of $i$-symmetric spaces with
orthogonal sum as addition modulo the so called {\it neutral spaces}. These groups are $4$-periodic:
$\W^{i}(\sheaf{T},\vee)\simeq\W^{i+4}(\sheaf{T},\vee)$.

\smallbreak

A {\it duality preserving functor} from $\sheaf{T}$ to another triangulated category with $\delta_{1}$-exact duality
$(\sheaf{T}_{1},\vee_{1},\varpi_{1})$ is a pair $(F,\eta)$, where $F:\sheaf{T}\too\sheaf{T}_{1}$ is an exact
functor and $\eta$ is a natural isomorphism $F\circ\vee\xrightarrow{\simeq}\vee_{1}\circ F$ satisfying
$$
\eta_{M^{\vee}}\circ F(\varpi_{M})\, =\, (\eta_{M})^{\vee_{1}}\circ\varpi_{1\, FM}
$$
and the equation $T_{1}^{-1}(\eta_{M})=(\delta_{1}\delta)\cdot\eta_{TM}$, where~$T_{1}$ denotes the translation
functor in~$\sheaf{T}_{1}$. The duality preserving functor~$(F,\eta)$ induces a homomorphism of triangular Witt groups,
see~\cite[Thm.\ 2.7]{Gi02}: If $(M,\varphi)$ is a $i$-symmetric space in $\sheaf{T}$ then
$$
(F,\eta)_{\ast}(M,\varphi)\, :=\; (F(M),\, (\delta\delta_{1})^{i}\cdot T_{1}^{i}(\eta_{M})\circ\varphi)
$$
is a $i$-symmetric space in~$\sheaf{T}_{1}$, which is neutral in $(\sheaf{T}_{1},\vee_{1},\varpi_{1})$
if $(M,\varphi)$ is neutral in the triangulated catgory with duality $(\sheaf{T},\vee,\varpi)$.

\smallbreak

Duality preserving functors can be composed, see~\cite[Sect.\ 1]{Gi03}.
Let $(G,\theta):(\sheaf{T}_{1},\vee_{1},\varpi_{1})\too (\sheaf{T}_{2},\vee_{2},\varpi_{2})$
be another duality preserving functor. The {\it composition} of~$(F,\eta)$ and $(G,\theta)$
is defined as follows:
\begin{equation}
\label{dualfunctorCompositionEq}
(G,\theta)\circ (F,\eta)\, :=\;\big(\, G\circ F\, ,\;\theta_{F}\circ G(\eta)\,\big)\, :\;
(\sheaf{T},\vee,\varpi)\,\too\, (\sheaf{T}_{2},\vee_{2},\varpi_{2})\, .
\end{equation}
We have then
$$
\big[\, (G,\theta)\circ (F,\eta)\,\big]_{\ast}(M,\varphi)\,\simeq\,
(G,\theta)_{\ast}\big(\, (F,\eta)_{\ast}(M,\varphi)\,\big)
$$
for all $i$-symmetric spaces~$(M,\varphi)$ in $(\sheaf{T},\vee,\varpi)$ and all $i\in\Z$.

\smallbreak

Another important definition is the following: Two duality preserving functors
$$
(F,\eta)\, ,\; (G,\theta)\, :\; (\sheaf{T},\vee,\varpi)\too (\sheaf{T}_{1},\vee_{1},\varpi_{1})
$$
are called {\it isometric} if there exists an isomorphism of functors $s:F\xrightarrow{\simeq}G$,
called {\it isometry}, which commutes with the respective translation functors and satisfies
$$
(s_{M})^{\vee_{1}}\circ\theta_{M}\circ s_{M^{\vee}}\, =\,\eta_{M}
$$
for all $M\in\sheaf{T}$. Then~$s_{M}$ is an isometry
$(F,\eta)_{\ast}(M,\varphi)\simeq (G,\theta)_{\ast}(M,\varphi)$.
\end{emptythm}

\begin{emptythm}
\label{DerWGrSubSect}
{\bf Derived Witt groups.}
The main example of a triangulated category with duality is the following.
Let~$(\sheaf{E},\vee,\varpi)$ be an exact category with duality. The the
derived functor of~$\vee$, which we denote (by some abuse of notation)
also by~$\vee$, is a duality on the bounded derived catgeory~$\mathrm{D}^{b}(\sheaf{E})$,
giving this triangulated category the structure of a triangulated category with duality. (The
isomorphism to the bidual is given in degree~$i$ by $\varpi_{M_{i}}:M_{i}\too M_{i}^{\vee\vee}$
for all $M_{\smb}\in\mathrm{D}^{b}(\sheaf{E})$.) The associated triangular Witt groups
are denote $\W^{i}(\sheaf{E},\vee)$, $i\in\Z$, and called the {\it derived Witt groups}
of $(\sheaf{E},\vee,\varpi)$.

\smallbreak

As usual in derived and coherent Witt theory we work with homological complexes.

\smallbreak

Note that by the main result of Balmer~\cite{Ba01a} we have an isomorphism
$$
\W_{\epsilon}(\sheaf{E},\vee)\,\xrightarrow{\;\simeq\;}\,\W^{1-\epsilon}(\sheaf{E},\vee)
$$
for all $\epsilon\in\{\pm 1\}$.
\end{emptythm}

\goodbreak
\section{Azumaya algebras with involutions and derived Witt groups}
\label{AzInvSect}\bigbreak

\begin{emptythm}
\label{NotationSubSect}
{\bf Notations and conventions.}
Let~$X$ be a noetherian scheme. We denote the structure sheaf of~$X$ by~$\OX$,
the local ring at~$x\in X$ by~$\OXx$, the maximal ideal of~$\OXx$ by~$\maxid_{x}$, and
set $k(x):=\OXx/\maxid_{x}$.

\smallbreak

Given an $\OX$-algebra~$\Az$ we use the following notations for categories of $\Az$-modules,
by which we mean -- if not otherwise said -- left $\Az$-modules: $\Mqc{\Az}$ the category of
quasi-coherent $\Az$-modules, and $\Mfg{\Az}$ the category of coherent $\Az$-modules.
We use also affine notations: If~$X=\Spec R$ we denote by~$A$ the global sections of~$\Az$
and write $\Mfg{A}$ and $\Mqc{A}$ instead of~$\Mfg{\Az}$ and $\Mqc{\Az}$, respectively.
\end{emptythm}

\begin{emptythm}
\label{AzAlgInvSubSect}
{\bf Azumaya algebras with involutions.}
By an {\it involution} of an $\OX$-algebra~$\Az$ we understand a $\OX$-linear homomorphism
$\inv:\Az\too\Az$ satisfying (i)~$\inv\circ\inv=\id_{\Az}$, and~(ii) $\inv_{U} (a\cdot b)=\inv_{U}(b)\cdot\inv_{U}(a)$
for all $a,b\in\Az (U)$ and all open $U\subseteq X$.

\smallbreak

Given an involution~$\inv$ on an $\OX$-algebra~$\Az$ we can turn a right $\Az$-module~$\sheaf{F}$
into a left one as follows: $a.x:=x\cdot\inv_{U}(a)$ for all $a\in\Az (U)$ and $x\in\sheaf{F}(U)$, $U\subseteq X$
open. We denoted this left $\Az$-module by $\ol{\sheaf{F}}$, or $\ol{\sheaf{F}}^{\inv}$, if we have to specify
the involution. Analogous we can turn a left $\Az$-module into a right one.

\smallbreak

Let~$R$ be a commutative ring and $(A,\inv)$ an $R$-algebra with involution.
We say that the pair $(A,\inv)$ is an {\it Azumaya algebra with involution} over~$R$ if
$A$ is a separable $R$-algebra, which is finitely generated and projective as $R$-module, and
the centre~$\ceZ (A)$ of~$A$ is either~$R$, in which case~$\inv$ is called of the {\it first kind},
or~$\ceZ (A)$ is a quadratic \'etale extension of~$R$ and~$R$ is the fix ring of~$\inv$, in which
case~$\inv$ is said to be of the {\it second kind}.

\smallbreak

Given a scheme~$X$ and an $\OX$-algebra~$\Az$ with involution~$\inv$ we say that
the pair $(\Az,\inv)$ is an {\it Azumaya algebra with involution of the first- or second kind over~$X$} if
it is locally an Azumaya algebra with involution of this kind.
\end{emptythm}

\begin{emptythm}
\label{HermWGrSubSect}
{\bf Derived hermitian Witt groups.}
Let~$(\Az,\inv)$ be an Azumaya algebra with involution (of first- or second kind) over the scheme~$X$,
and $\Pcat{\Az}$ the full subcategory of $\Mfg{\Az}$ consisting of coherent $\Az$-modules, which are
locally free as $\OX$-modules. The contravariant functor
$$
\dual{\Az,\inv}\, :\;\sheaf{F}\,\longmapsto\,\ol{\HOM_{\Az}(\sheaf{F},\Az)}\, =\,
\ol{\HOM_{\Az}(\sheaf{F},\Az)}^{\inv}
$$
is a duality on~$\Pcat{\Az}$ making this an exact category with duality. Associated with this data we have
the 'classical' (skew-)hermitian Witt groups~$\W_{\epsilon}(\Az,\inv)$, $\epsilon=\pm 1$, and the derived
Witt groups, denoted $\W^{i}(\Az,\inv)$, $i\in\Z$, called the {\it derived hermitian Witt groups} of $(\Az,\inv)$.

\medbreak

Let $f:Y\too X$ be a morphism of schemes. The natural isomorphism of left $f^{\ast}\Az$-modules
$f^{\ast}\ol{\HOM_{\Az}(\sheaf{F},\Az)}\xrightarrow{\simeq}\ol{\HOM_{f^{\ast}\Az}(f^{\ast}\sheaf{F},f^{\ast}\Az)}$
makes the pull-back functor $f^{\ast}:\DP{\Az}\too\DP{f^{\ast}\Az}$ duality preserving and therefore we have
a homomorphism $f^{\ast}:\W^{i}(\Az,\inv)\too\W^{i}(f^{\ast}(\Az,\inv))$ for all~$i\in\Z$.

\medbreak

There are also Witt groups with support. Recall first that the support $\supp\sheaf{F}_{\smb}$
of a complex~$\sheaf{F}_{\smb}$ of $\OX$-modules is the set of all $x\in X$ with
$\HM_{i}(\sheaf{F}_{\smb\, x})\not= 0$ for at least one~$i\in\Z$. Here $\HM_{i}(\sheaf{F}_{\smb\, x})$
denotes the $i$th homology group of the at~$x\in X$ localized complex~$\sheaf{F}_{\smb\, x}$.

\smallbreak

Let~$\DPS{\Az}{Z}$ be the full triangulated subcategory of $\DP{\Az}$ consisting of complexes~$\sheaf{F}_{\smb}$
with support in the closed subscheme  $Z\subseteq X$. The restriction of $\dual{\Az,\inv}$ to this category
makes $\DPS{\Az}{Z}$ a triangulated category with duality. Its associated triangular Witt groups are the so called
{\it derived hermitian Witt groups} of~$(\Az,\inv)$ {\it with support} in~$Z$, denoted $\W^{i}_{Z}(\Az,\inv)$, $i\in\Z$.

\smallbreak

As usual if~$X=\Spec R$ and~$Z$ is defined by an ideal~$\ideal$ we set $A:=\Gamma (X,\Az)$ and use
affine notations $\DPS{A}{\ideal}$ instead of $\DPS{\Az}{Z}$, and $\W^{i}_{\ideal}(A,\inv)$ instead of
$\W^{i}_{Z}(\Az,\inv)$.
\end{emptythm}

\goodbreak
\section{Derived categories and dualizing complexes}
\label{DualComplSect}\bigbreak

\begin{emptythm}
\label{dualComplNotationSubSect}
{\bf Some (derived) categories of modules.}
Let~$R$ be a commutative noetherian ring, $\ideal\subset R$~an ideal, and~$A$
an $R$-algebra. We denote by $\MqcS{A}{\ideal}$ the category of $A$-modules with support in the
closed subscheme $\Spec (R/\ideal)\subseteq\Spec R$, and set
$\MfgS{A}{\ideal}=\Mfg{A}\cap\MqcS{A}{\ideal}$.

\smallbreak

We denote further by~$\Dfg{A}$ the full subcategory of the bounded derived category~$\DMqc{A}$ consisting
of complexes with coherent homology, and by $\DcS{A}{\ideal}$ the full subcategory of $\Dfg{A}$
consisting of complexes with support in~$\Spec (R/\ideal)\subseteq\Spec R$.

\smallbreak

For~$p\in\Z$ we denote by~$\DfgF{A}{p}$ the full subcategory of~$\Dfg{A}$
consisting of complexes~$M_{\smb}$ with $\codimension_{\Spec R}\supp M_{\smb}\geq p$,
and set $\DcSF{A}{\ideal}{p}:=\DcS{A}{\ideal}\cap\DfgF{A}{p}$.

\smallbreak

We have the following well known equivalence of derived categories, which
follows from a result due to Grotendieck, which is proven in Verdier's thesis~\cite[pp 169--170]{Ve96},
see~\cite[Sect.\ 4.2]{Gi13} for details.
\end{emptythm}

\begin{emptythm}
\label{derCatEqLem}
{\bf Lemma.}
{\it
Let~$R$ be a commutative noetherian ring with ideal~$\ideal$ and~$A$ a coherent
$R$-algebra. Then the natural functor $\mathrm{D}^{b}(\MfgS{A}{\ideal})\too\DcS{A}{\ideal}$
is an equivalence.
}
\end{emptythm}

\begin{emptythm}
\label{signConventionsSubSect}
{\bf Sign conventions.}
Let~$R$ and~$A$ be as above, and
$I_{\smb}\in\Dfg{R}$ a bounded complex of injective $R$-modules.

\smallbreak

If~$M$ is a left $A$-module and $N$~an arbitrary $R$-module
then $\Hom_{R}(M,N)$ becomes a right $A$-module by setting
$(f\cdot a)(m):=f(am)$ for $f\in\Hom_{R}(M,N)$, $a\in A$, and $m\in M$.
Analogous if~$M$ is a right $A$-module then $\Hom_{R}(M,N)$ has a
left $A$-module structure: $(a\cdot f)(m):=f(ma)$.

\smallbreak

For~$M_{\smb}\in\Dfg{R}$ the complex $\Hom_{R}(M_{\smb},I_{\smb})$ is given in degree~$l$ by
$$
\Hom_{R}(M_{\smb},I_{\smb})_{l}\, =\;\bigoplus\limits_{r\in\Z}\Hom_{R}(M_{-l-r},I_{-r})\, ,
$$
and the $r$-component of the differential $\Hom_{R}(M_{\smb},I_{\smb})_{l}\too\Hom_{R}(M_{\smb},I_{\smb})_{l-1}$
maps $g\in\Hom_{R}(M_{-l-r},I_{-r})$ onto $g\circ d^{M}_{-l-r+1}+(-1)^{l+1}d^{I}_{-r}\circ g$,
where $d^{M}_{\smb}$ and $d^{I}_{\smb}$ denote the differentials of~$M_{\smb}$ and $I_{\smb}$, respectively.

\smallbreak

The natural homomorphism $\varpi^{I}_{M}\, :\; M_{\smb}\too\Hom_{R}(\Hom_{R}(M_{\smb},I_{\smb}),I_{\smb})$
is defined as follows: The $(r,s)$-component of
$$
(\varpi^{I}_{M})_{l}\, :\; M_{l}\,\too\,\bigoplus\limits_{r,s\in\Z}\Hom_{R}(\Hom_{R}(M_{l+s-r},I_{-r}),I_{-s})
$$
is~$0$ if $r\not= s$, and otherwise equal $(-1)^{\frac{s(s+1)}{2}}$
times the evaluation map $M_{l}\too\Hom_{R}(\Hom_{R}(M_{l},I_{-s}),I_{-s})$.
Note that the evaluation map is (left-) $A$-linear if~$M_{\smb}$ is a complex
of (left-) $A$-modules.
\end{emptythm}

\begin{emptythm}
\label{DualComplDef}
{\bf Definition.}
The complex~$I_{\smb}$ is called a {\it dualizing complex} of~$R$ if $\varpi^{I}_{M}$
is an isomorphism in~$\Dfg{R}$ for all $M_{\smb}\in\Dfg{R}$.
\end{emptythm}

\begin{emptythm}
\label{dualComplExpl}
{\bf Examples.}
\begin{itemize}
\item[(i)]
Let~$R$ be a Gorenstein ring of finite Krull dimension. Then a finite injective resolution
$I_{0}\too I_{-1}\too\ldots\too I_{-\dim R}$ of an invertible $R$-module considered as
an element of~$\Dfg{R}$ with $I_{l}$ in degree~$l$ is a dualizing complex of~$R$.

\smallbreak

\item[(ii)]
Let $\alpha:R\too S$ be a finite morphism of noetherian rings and assume~$R$ has a dualizing complex~$I_{\smb}$.
Then
$$
\alpha^{\natural}(I_{\smb})\, :=\;\Hom_{R}(S,I_{\smb})\;\in\,\Dfg{S}
$$
is a dualizing complex of~$S$. This can be seen as follows. The morphism of complexes
$\Hom_{R}(S,I_{\smb})\too I_{\smb}$ given in degree~$l$ by $h\mapsto (-1)^{l}h(1)$, induces a natural
isomorphism in~$\Dfg{R}$
$$
\vartheta_{M}\, :\;\Hom_{S}(M_{\smb},\Hom_{R}(S,I_{\smb}))\,\xrightarrow{\;\simeq\;}\,
\Hom_{R}(M_{\smb},I_{\smb})\, ,
$$
which implies that
$\varpi^{\alpha^{\natural}(I)}_{M}\, :\; M_{\smb}\too\Hom_{S}(\Hom_{S}(M_{\smb},\alpha^{\natural}(I_{\smb})),\alpha^{\natural}(I_{\smb}))$
is an isomorphism in~$\Dfg{S}$ for all $M_{\smb}\in\Dfg{S}$.

\smallbreak

In case~$I_{\smb}:I_{0}\too I_{-1}\too\ldots$ is a finite injective resolution of the $R$-module~$R$,
$t\in R$~a non unit and non zero divisor, and $\alpha:R\too S:=R/Rt$ the quotient morphism then
$\alpha^{\natural}(I_{\smb})$ is an injective resolution of $R/Rt$ living in degrees $-1,\ldots ,1-\dim R$.
This can be seen using the long exact $\Ext_{R}(\, -\, ,R)$-sequence associated with the short exact sequence
$0\too R\xrightarrow{\cdot t} R\too R/Rt\too 0$. (In this case~$R$ and $S=R/Rt$ are necessarily Gorenstein
rings of finite Krull dimension.)
\end{itemize}
\end{emptythm}

\begin{emptythm}
\label{codfctSubSect}
{\bf The codimension function of a dualizing complex.}
Let~$I_{\smb}$ be a  dualizing complex of the ring~$R$. Then given $P\in\Spec R$ there exists precisely
one integer $l=l(P)$, such that $\HM_{d}(\Hom_{R_{P}}(k(P),I_{\smb\, P}))=0$ for all $d\not= l$
and $\HM_{l}(\Hom_{R_{P}}(k(P),I_{\smb\, P}))\simeq k(P)$, where~$k(P)$ is the residue field of~$P$, see~\cite[Chap.\ V, Prop.\ 3.4]{RD}.
This integer will be denoted $-\cod{I}(P)$ (for the minus sign, note that we use homological complexes). We get a function
$$
\cod{I}\, :\;\Spec R\,\too\,\Z\, ,\; P\,\longmapsto\,\cod{I}(P)\, ,
$$
the {\it codimension function} of the dualizing complex~$I_{\smb}$. For instance, if~$R$ is a Gorenstein
ring of finite Krull dimension and~$I_{\smb}:I_{0}\too I_{-1}\too\ldots$ is an injective resolution of a rank
one projective $R$-module with $I_{0}$ in degree~$0$ then $\cod{I}(P)=\dim R_{P}$ for
all $P\in\Spec R$, see \eg~\cite[Sect.\ 3.3]{CMR}.
\end{emptythm}

\begin{emptythm}
\label{dualComplLem}
{\bf Lemma.}
{\it
Let~$R$ be a Gorenstein domain of finite Krull dimension and $I_{\smb}:I_{0}\too I_{-1}\too\ldots\too I_{-\dim R}$
a injective resolution of an invertible $R$-module~$L$ considered as an element of~$\Dfg{R}$ with~$I_{0}$ in
degree~$0$. If $\alpha:R\too S$ is a finite morphism of rings with $\dim S=\dim R$ then
$\cod{\alpha^{\natural}(I)}(Q)=\dim S_{Q}$ for all $Q\in\Spec S$.
}

\begin{proof}
Since~$R$ and~$S$ have the same Krull dimension, and since $\dim R/J<\dim R$ for all non zero
ideals~$J$ of~$R$ the morphism~$\alpha$ is injective. Let~$Q$ be a prime ideal of~$S$ and
$P=\alpha^{-1}(Q)$. Replacing~$R$ by~$R_{P}$ we can assume that~$R$ is a local ring with
maximal ideal~$P$, and so~$Q$ is a maximal ideal of~$S$ as well.

\smallbreak

We have then $\dim S_{Q}=\dim R=\cod{I}(P)$, and $S/Q\simeq (R/P)^{\oplus\, m}$
for some integer~$m\geq 1$. As seen above, see Example~\ref{dualComplExpl}~(ii), we have
quasi-isomorphisms of complexes of $R$-modules
$$
\Hom_{S}(S/Q,\alpha^{\natural}(I_{\smb}))\,\simeq\,\Hom_{R}(S/Q,I_{\smb})\,\simeq\,
\Hom_{R}(R/P,I_{\smb})^{\oplus\, m}\, .
$$
Since $\HM_{l}(\Hom_{R}(R/P,I_{\smb}))\not=0$ if and only if $l=-\dim R$
it follows that $\cod{\alpha^{\natural}(I)}(Q)=\dim R=\dim S_{Q}$. We are done.
\end{proof}
\end{emptythm}

\goodbreak
\section{Coherent hermitian Witt groups}
\label{CohWGrSect}\bigbreak

\begin{emptythm}
\label{CohHermWGrSubSect}
{\bf Definition of coherent hermitian Witt groups.}
We refer to~\cite{Gi07,Gi13} for details and more information.

\smallbreak

Let~$R$ be a commutative ring with dualizing complex~$I_{\smb}\in\Dfg{R}$, and $(A,\inv)$ a
coherent $R$-algebra with involution~$\inv$. The derived functor
$$
\dual{A,\inv}_{I}\, :\;M_{\smb}\,\longmapsto\,\ol{\Hom_{R}(M_{\smb},I_{\smb})}\, =\,
\ol{\Hom_{R}(M_{\smb},I_{\smb})}^{\inv}
$$
is a duality on~$\Dfg{A}$ making this a triangulated category with duality. The isomorphism
to the bidual is given by~$\varpi^{I}$, which is a quasi-isomorphism of complexes of $A$-modules.

\smallbreak

The associated triangular Witt groups are denoted $\cw{A,\inv,I_{\smb}}{i}$, $i\in\Z$, and called {\it coherent
hermitian Witt groups} of~$(A,\inv)$. There are also Witt groups with support. If~$\ideal\subseteq R$
is an ideal then~$\DcS{A}{\ideal}$ is also a triangulated category with duality $\dual{A,\inv}_{I}$. The associated
triangular Witt groups $\cws{A,\inv,I_{\smb}}{\ideal}{i}$, $i\in\Z$, are called {\it coherent hermitian Witt groups}
of~$(A,\inv)$ {\it with support} in the ideal~$\ideal$.
\end{emptythm}

\begin{emptythm}
\label{der-cohWGrSubSect}
{\bf Derived and coherent Witt groups.}
Assume now that~$(A,\inv)$ is an Azumaya algebra with involution of first- or second kind, and
that~$R$ is a regular ring of finite Krull dimension. Then a (finite) injective resolution
$I_{\smb}: I_{0}\too I_{-1}\too\ldots\too I_{-\dim R}$ of~$R$ considered as an element of $\Dfg{R}$ with
$I_{-i}$ in degree~$-i$ is a dualizing complex of~$R$. Under these assumptions
the natural functor $\DP{A}\too\Dfg{A}$ is an equivalence, which becomes duality preserving via
$$
\Hom_{A}(\, -\, ,A)\,\xrightarrow{\;\simeq\;}\,\Hom_{R}(\, -\, ,R)\,\xrightarrow{\;\simeq\;}\,
\Hom_{R}(\, -\, ,I_{\smb})\, ,
$$
where the isomorphism of functors on the left hand side is induced by the reduced trace composed with
the standard trace of a quadratic \'etale extension if~$\inv$ is of the second kind, and the other one by a
quasi-isomorphism $R\xrightarrow{\simeq}I_{\smb}$. We refer to~\cite[App.]{Gi09} for proofs and details.

\smallbreak

In particular, we have then isomorphisms $\W^{i}(A,\inv)\xrightarrow{\simeq}\cw{A,\inv,I_{\smb}}{i}$ for
all~$i\in\Z$, and analogous for the Witt groups with support.
\end{emptythm}

\begin{emptythm}
\label{transferSubSect}
{\bf The transfer map.}
Let $\alpha:R\too S$ be a finite homomorphism of commutative noetherian rings,
and $I_{\smb}$ a dualizing complex of~$R$. Let further~$(A,\inv)$ be a coherent $R$-algebra with
involutions and $(B,\binv)$ a coherent $S$-algebra with involution. Assume that there exists
a $R$-algebra homomorphism $\xi:A\too B$, which is compatible with the involutions, \ie we have
$\xi\circ\inv=\binv\circ\xi$. We indicate this situation by writing
\begin{equation}
\label{pairMorphismEq}
(\alpha,\xi)\, :\;\big(\, R,(A,\inv)\,\big)\,\too\,\big(\,S,(B,\binv)\,\big)\, .
\end{equation}

\smallbreak

\noindent
Then $\alpha^{\natural}(I_{\smb}):=\Hom_{R}(S,I_{\smb})$ is a dualizing complex of~$S$ and
the morphism of functors~$\eta$ introduced in Example~\ref{dualComplExpl}~(ii) induces a
natural quasi-isomorphism of complexes of $R$-modules:
$$
\vartheta^{\xi}_{M}:
\alpha_{\ast}\big(\dual{B,\binv}_{\alpha^{\natural}(I)}(M_{\smb})\big)=
\ol{\Hom_{S}(M_{\smb},\alpha^{\natural}(I_{\smb}))}\too\ol{\Hom_{R}(M_{\smb},I_{\smb})}
=\dual{A,\inv}_{I}\big(\alpha_{\ast}(M_{\smb})\big)\, .
$$
Since~$\xi:A\too B$ is $R$-linear and compatible with the involutions, the quasi-isomorphism
of complexes $\vartheta^{\xi}_{M}$ is a morphism of complexes of $A$-modules, \ie
an isomorphism in~$\Dfg{A}$ for all $M_{\smb}\in\Dfg{B}$.

\smallbreak

A straightforward verification shows that~$\vartheta^{\xi}$ is a duality transformation for the push-forward
$\alpha_{\ast}\, :\;\Dfg{B}\,\too\,\Dfg{A}$. Therefore given a $i$-symmetric space $(M_{\smb},\varphi)$ in
$(\Dfg{B},\dual{B,\binv}_{\alpha^{\natural}(I)})$ then
$$
\tr_{(\alpha,\xi)}(M_{\smb},\varphi)\, :=\; (\alpha_{\ast},\vartheta^{\xi})_{\ast}(M_{\smb},\varphi)
$$
is a $i$-symmetric space in $(\Dfg{A},\dual{A,\inv}_{I})$.

\medbreak

Let now $(\beta,\xi_{1}):\big(S,(B,\binv)\big)\too\big(S_{1},(B_{1},\binv_{1})\big)$
be another morphism, where~$\beta$ is a finite morphism and $(B_{1},\binv_{1})$
is a coherent $S_{1}$-algebra with involution. Then we define
$$
(\beta,\xi_{1})\circ (\alpha,\xi)\, :=\; (\beta\circ\alpha,\xi_{1}\circ\xi)\, :\;
\big(\, R,(A,\inv)\,\big)\,\too\,\big(\, S_{1},(B_{1},\binv_{1})\,\big)\, .
$$
Identifying $(\beta\circ\alpha)^{\natural}(I_{\smb})=\alpha^{\natural}\big(\beta^{\natural}(I_{\smb})\big)$
the identity $(\beta\circ\alpha)_{\ast}=\alpha_{\ast}\circ\beta_{\ast}$ is an
isometry of duality preserving functors
$$
\big((\beta\circ\alpha)_{\ast},\vartheta^{\xi_{1}\circ\xi}\big)\,\xrightarrow{\;\simeq\;}\,
(\alpha_{\ast},\vartheta^{\xi})\circ (\beta_{\ast},\vartheta^{\xi_{1}})\, ,
$$
and so we have an isometry
\begin{equation}
\label{TrIsometryEq}
\tr_{(\beta\circ\alpha\, ,\;\xi_{1}\circ\xi)}(M_{\smb},\varphi)\,\simeq\
\tr_{(\alpha,\xi)}\big(\tr_{(\beta,\xi_{1})}(M_{\smb},\varphi)\big)
\end{equation}
in $(\Dfg{A},\dual{A,\inv}_{I})$ for all $i$-symmetric spaces $(M_{\smb},\varphi)$ in the triangulated
category with duality $(\Dfg{B_{1}},\dual{B_{1},\binv_{1}}_{(\beta\circ\alpha)^{\natural}(I)})$.

\medbreak

\noindent
{\bf Example.}
Let as above $\alpha:R\too S$ be a finite
homomorphism of noetherian rings, where~$R$ has a dualizing complex~$I_{\smb}$, and $(A,\inv)$
a coherent $R$-algebra with involution. Then $(B,\binv):=S\otimes_{R}(A,\inv)$ is a coherent
$S$-algebra with involution and we have a natural morphism of $R$-algebras
$$
\xi\, :\; A\,\too\, B\, ,\; a\,\longmapsto\, 1\otimes a\, ,
$$
which is compatible with the involutions. We get a duality
preserving functor $(\alpha_{\ast},\vartheta^{\xi})$ and a transfer map $\tr_{(\alpha,\xi)}$, which
we denote also~$\tr_{\alpha}$ only as there is a canonical choice for the duality transformation.
\end{emptythm}

\begin{emptythm}
\label{devSubSect}
{\bf D\'evissage.}
Let~$R$ be a Gorenstein ring of finite Krull dimension, $I_{\smb}:I_{0}\too I_{-1}\too\ldots\in\Dfg{R}$
a finite injective resolution of the $R$-module~$R$ living in the indicated degrees and~$(A,\inv)$ a
coherent $R$-algebra with involution. Let further~$t\in R$ be a non unit and non zero divisor, and
$\pi:R\too R/Rt$ the quotient morphism. We set $(A',\inv'):=R/Rt\otimes_{R}(A,\inv)$.

\medbreak

Then it is shown in~\cite[Sect.\ 5]{Gi07} that mapping an $i$-symmetric space
$(M_{\smb},\varphi)$ in $(\DfgF{A'}{p-1},\dual{A',\inv'}_{\pi^{\natural}{I}})$ onto
$\tr_{\pi}(M_{\smb},\varphi)$ induces an isomorphism
$$
\W^{i}(\DfgF{A'}{p},\dual{A',\inv'}_{\pi^{\natural}{I}})\,\xrightarrow{\;\simeq\;}\,
\W^{i}(\DcSF{A}{Rt}{p+1},\dual{A,\inv}_{I})
$$
for all~$i\in\Z$ and integers~$p\geq 0$. (Note here that in~\cite{Gi07} the filtrations on the
bounded derived categories are defined using the codimension functions associated with
the respective dualizing complexes. But we have $\cod{\pi^{\natural}(I_{\smb})}(P/Rt)=\hgt P$
for all prime ideals~$P\supseteq Rt$, see Example~\ref{dualComplExpl}~(ii).)
\end{emptythm}

\begin{emptythm}
\label{PullBackSubSect}
{\bf Pull-backs.}
Let $\alpha:R\too S$ be a flat morphism of commutative noetherian rings with dualizing complexes
$I_{\smb}\in\Dfg{R}$ and $J_{\smb}\in\Dfg{S}$, respectively, and~$(A,\inv)$ be
a coherent $R$-algebra with involution. We set $(B,\binv):=S\otimes_{R}(A,\inv)$.

\smallbreak

Assume that there exists a quasi-isomorphism of complexes of $S$-modules
$$
\rho\, :\; S\otimes_{R}I_{\smb}\,\too\, J_{\smb}\, .
$$
This quasi-isomorphism induces a natural (in~$M_{\smb}$) quasi-isomorphism $c_{\rho\, M}$:
$$
S\otimes_{R}\dual{A,\inv}_{I}(M_{\smb})\, =\,
S\otimes_{R}\ol{\Hom_{R}(M_{\smb},I_{\smb})}\,\xrightarrow{\simeq}\,
\ol{\Hom_{S}(S\otimes_{R}M_{\smb},S\otimes_{R}I_{\smb})}\qquad\qquad
$$
$$
\qquad\qquad\qquad\qquad
\xrightarrow{\;\;\;\ol{\Hom_{S}(S\otimes M,\rho)}\;\;\;}\,\ol{\Hom_{S}(S\otimes_{R}M_{\smb},J_{\smb})}\,
=\,\dual{B,\binv}_{J}(S\otimes_{R}M_{\smb})
$$
for all $M_{\smb}\in\Dfg{A}$, which is a duality transformation for the pull-back~$\alpha^{\ast}$.
We get a duality preserving functor
$$
(\alpha^{\ast},c_{\rho})\, :\; (\Dfg{A},\dual{A,\inv}_{I})\,\too\, (\Dfg{B},\dual{B,\binv}_{J})\, .
$$

\medbreak

\noindent
{\bf Example.}
Let $R$ be a commutative noetherian ring, $I_{\smb}$ a dualizing complex
of~$R$, $(A,\inv)$ a coherent $R$-algebra with involution, and $\alpha:R\too S$ an open
immersion, \ie the induced morphism of affine schemes $\Spec S\too\Spec R$ is an open immersion,
or a localization at some multiplicative closed subset of~$R$.

\smallbreak

Then $\alpha^{\ast}(I_{\smb})=S\otimes_{R}I_{\smb}$ is a dualizing complex
of~$S$, and we have a canonical pull-back
$$
(\alpha^{\ast},c_{\,\id_{S\otimes I}})\, :\;
(\Dfg{A},\dual{A,\inv}_{I})\,\too\,(\Dfg{S\otimes_{R}A},\dual{S\otimes (A,\inv)}_{S\otimes I})\, ,
$$
which we denote by~$\alpha^{\ast}$ only.

\medbreak

The proof of the following result, which generalizes~\cite[Lem.\ 3.5]{Gi13}, is straightforward.
\end{emptythm}

\begin{emptythm}
\label{Pushforward-PullbackLem}
{\bf Lemma.}
{\it
Let $\alpha:R\too S$ be a flat morphism of commutative noetherian rings with dualizing complexes
$I_{\smb}\in\Dfg{R}$ and $J_{\smb}\in\Dfg{S}$, respectively, $\pi_{R}:R\too\bar{R}$ a finite morphism, and
$$
\xymatrix{
\bar{R} \ar[r]^-{\bar{\alpha}}  & \bar{S}
\\
R \ar[r]_-{\alpha} \ar[u]^-{\pi_{R}} & S \ar[u]_-{\pi_{S}} \rlap{\, ,}
}
$$
the corresponding cartesian square, \ie $\bar{S}=\bar{R}\otimes_{R}S$.
Let further $(A,\inv)$ be an Azumaya algebra with involution over~$R$.

\smallbreak

If there exists a quasi-isomorphism of complexes of $S$-modules $\rho: S\otimes_{R}I_{\smb}\too J_{\smb}$ then

\smallbreak

\begin{itemize}
\item[(i)]
the morphism of complexes of $\bar{S}$-modules $\bar{\rho}$:
$$
\bar{S}\otimes_{\bar{R}}\pi_{R}^{\natural}(I_{\smb})\, =\, S\otimes_{R}\bar{R}\otimes_{\bar{R}}\Hom_{R}(\bar{R},I_{\smb})
\,\xrightarrow{\;\simeq\;}\Hom_{S}(S\otimes_{R}\bar{R},S\otimes_{R}I_{\smb})
$$
$$
\qquad\qquad\xrightarrow{\;\;\;\Hom_{S}(\bar{S},\rho)\;\;\;}\,\Hom_{S}(\bar{S},J_{\smb})\, =\,\pi_{S}^{\natural}(J_{\smb})
$$
is a quasi-isomorphism, and

\smallbreak

\item[(ii)]
the natural isomorphism of functors
$\alpha^{\ast}\circ\pi_{R\,\ast}\xrightarrow{\simeq}\pi_{S\,\ast}\circ\bar{\alpha}^{\ast}$
is an isometry of duality preserving functors
$$
(\alpha^{\ast}, c_{\rho})\circ (\pi_{R\,\ast},\vartheta^{\xi_{R}})\,\xrightarrow{\;\simeq\:}\,
(\pi_{S\,\ast},\vartheta^{\xi_{S}})\circ (\bar{\alpha}^{\ast},c_{\bar{\rho}})\, ,
$$
where $\xi_{R}:A\too\bar{R}\otimes_{R}A$ and $\xi_{S}:S\otimes_{R}A\too\bar{S}\otimes_{R}A$
are the natural morphisms of $R$- respectively $S$-algebras.
\end{itemize}
}
\end{emptythm}

\begin{emptythm}
\label{zeroThmSubSect}
{\bf The zero theorem.}
Let
$$
\xymatrix{
 & \tilde{R}
 \\
S \ar[ru]^-{\beta} & R \ar[l]^-{u} \ar[u]_-{\gamma}
}
$$
be a commutative diagram of Gorenstein rings of finite Krull dimension,
where $u$ is flat, $\gamma$ a localization morphism, \ie~$\tilde{R}=U^{-1}R$
for some multiplicative closed subset~$U$ of~$R$, and~$\beta$ is a surjective morphism
with kernel generated by a non zero divisor~$y\in S$. Let further $(A,\inv)$ be an
Azumaya algebra with involution of the first- or second kind over~$R$,
$(B,\binv):=u^{\ast}(A,\inv)$, and $R\too I_{0}\too I_{-1}\too\ldots$ and
$S\too J_{0}\too J_{-1}\too\ldots$ finite injective resolutions of~$R$
and~$S$, respectively, living in the indicated degrees in~$\Dfg{R}$
and~$\Dfg{S}$, respectively.

\smallbreak

With this notation we have the (so called) zero theorem:

\medbreak

\noindent
{\bf Theorem (\cite[Thm.\ 6.3]{Gi13}).}
{\it
Given a $i$-symmetric space $(M_{\smb},\varphi)$ in the triangulated category with duality
$(\DfgF{A}{p},\dual{A,\inv}_{T^{-1}I})$, $p\geq 0$, then
the transfer $\tr_{\beta}\big(\,\gamma^{\ast}(M_{\smb},\varphi)\,\big)$
is neutral in $(\DfgF{B}{p},\dual{B,\binv}_{J})$.
}
\end{emptythm}

\goodbreak
\section{A technical lemma}
\label{techLemSect}\bigbreak

\begin{emptythm}
\label{maintechLem1SubSect}
Let~$R$ be a commutative noetherian ring with dualizing complex~$I_{\smb}$ and $U\subset R$ a
multiplicative closed subset. We denote $\alpha:R\too U^{-1}R$ the localization morphism. Let further
$\beta: U^{-1}R\too S$ be a morphism, such that $\beta\circ\alpha:R\too S$ is onto.
Set $\Jideal:=\Ker (\beta\circ\alpha)$. Then we have an isomorphism
$R/\Jideal\xrightarrow{\simeq}U^{-1}(R/\Jideal)\simeq U^{-1}R/U^{-1}\Jideal\simeq S$.

\smallbreak

Let further~$A$ be a coherent $R$-algebra.
\end{emptythm}

\begin{emptythm}
\label{catEqLem}
{\bf Lemma.}
{\it
Let $M\in\MqcS{A}{\Jideal}$ and $N\in\MfgS{U^{-1}A}{U^{-1}\Jideal}$. Then:

\smallbreak

\begin{itemize}
\item[(i)]
$M$ is a $U^{-1}A$-module and the localization
homomorphism $\iota^{M}:M\too U^{-1}M$ is an isomorphism of $U^{-1}A$-modules.

\smallbreak

\item[(ii)]
$N$ is finitely generated as $A$-module.

\smallbreak

\item[(iii)]
The pull-back $\alpha^{\ast}:\MfgS{A}{\Jideal}\too\MfgS{U^{-1}}{U^{-1}\Jideal}$ is
an equivalence with inverse the push-forward~$\alpha_{\ast}$.
\end{itemize}
}

\begin{proof}
(iii)~is a consequence of~(i) and~(ii). To prove~(i) we observe that every~$u\in U$ is
invertible modulo~$\Jideal$ and so also module~$\Jideal^{m}$ for all $m\geq 1$. In fact,
if $ur+x=1$ for some $r\in R$ and~$x\in\Jideal$ then $1=(ur+x)^{m}=u\cdot s+x^{m}$ for some~$s\in R$ by
the binomial formula. Hence given a $A$-module~$M$ with support in $\Spec R/\Jideal$
then~$M$ is a $U^{-1}A$-module and the natural homomorphism $\iota^{M}:M\too U^{-1}M$ is an
isomorphism of $U^{-1}A$-modules.

\smallbreak

Finally, we prove~(ii). By assumption there exists an integer~$l\geq 1$, such
that $\Jideal^{l}N=0$. If~$l=1$ then~$N$ is a finitely generated $R/\Jideal$-module and so also
finitely generated as $A$-module. If~$l\geq 2$ we conclude by induction using the exact
sequence $0\too\Jideal^{l-1}N\too N\too N/\Jideal^{l-1}N\too 0$.
\end{proof}

\smallbreak

By Lemma~\ref{derCatEqLem} this has the following implication.

\smallbreak

\noindent
{\bf Corollary.}
{\it
The pull-back
$$
\alpha^{\ast}\, :\;\DcS{A}{\Jideal}\,\too\,\DcS{U^{-1}A}{U^{-1}\Jideal}
$$
is an equivalence with inverse the push-forward $\alpha_{\ast}$. 
}
\end{emptythm}

\begin{emptythm}
\label{invDualityTrafoSubSect}
Assume now that~$(A,\inv)$ is an Azumaya algebra with involution over~$R$.

\smallbreak

The equivalence~$\alpha^{\ast}$ is duality preserving
with duality transformation the natural isomorphism
$$
c_{\id_{U^{-1}I}}\, :\; U^{-1}\ol{\Hom_{R}(\, -\, ,I_{\smb})}\,\xrightarrow{\;\simeq\;}\,
\ol{\Hom_{U^{-1}R}(U^{-1}\, -\, ,U^{-1}I_{\smb})}\, ,
$$
see  the example in~\ref{PullBackSubSect}.

\medbreak

Therefore by Balmer and Walter~\cite[Lem.\ 4.3~(d)]{BaWa02} the
inverse equivalence~$\alpha_{\ast}$ is duality preserving as well. A duality transformation for~$\alpha_{\ast}$
$$
\theta\, :\;\ol{\Hom_{U^{-1}R}(\, -\, ,U^{-1}I_{\smb})}\,\xrightarrow{\;\simeq\;}\,
\ol{\Hom_{R}(\, -\, ,I_{\smb})}\, ,
$$
is defined as the inverse of the isomorphism of complexes of $A$-modules
\begin{equation}
\label{ThetaDefEq}
\eta_{N}\, :\;\ol{\Hom_{R}(N_{\smb},I_{\smb})}\,\xrightarrow{\;\iota^{\Hom_{R}(N,I)}\;}\,
U^{-1}\ol{\Hom_{R}(N_{\smb},I_{\smb})}\,\xrightarrow{\;c_{\id_{U^{-1}I}}\;}\,
\qquad\qquad\qquad\qquad
\end{equation}
$$
\qquad\qquad\ol{\Hom_{U^{-1}R}(U^{-1}N_{\smb}, U^{-1}I_{\smb})}\,\xrightarrow{\;\Hom_{U^{-1}R}(\iota^{N},U^{-1}I)\;}\,
\ol{\Hom_{U^{-1}R}(N_{\smb},U^{-1}I_{\smb})}\, .
$$
for all $N_{\smb}\in\DcS{U^{-1}A}{U^{-1}\Jideal}$. Here~$A$ acts on the $U^{-1}A$-module~$N_{i}$, $i\in\Z$,
via the natural homomorphism of $R$-algebras $\xi:A\too U^{-1}A$, and the first and last morphism of complexes
is an isomorphism by Lemma~\ref{catEqLem}~(i) above.

\smallbreak

We observe that the quasi-isomorphism $\eta_{N}$ has in degree~$l$ the $r$-component
$$
\ol{Hom_{R}(N_{-l-r},I_{-r})}\,\too\,\ol{\Hom_{U^{-1}R}(N_{-l-r},U^{-1}I_{-r})}\, ,\;
h\,\longmapsto\,\iota^{I_{-r}}\circ h\, .
$$
\end{emptythm}

\begin{emptythm}
\label{maintechLem2SubSect}
We have the two dualizing complexes $(\beta\circ\alpha)^{\natural}(I_{\smb})$ and
$\beta^{\natural}(\alpha^{\ast})=\beta^{\natural}(U^{-1}I_{\smb})$ on~$S$. As seen
above the natural $S$-linear morphism
$$
(\beta\circ\alpha)^{\natural}(I_{\smb})\, =\,\Hom_{R}(S,I_{\smb})\,\too\,\Hom_{U^{-1}R}(S,U^{-1}I_{\smb})
\, =\,\beta^{\natural}(U^{-1}I_{\smb})\, ,
$$
given by $h\mapsto\iota^{I_{l}}\circ h$ in degree~$l$, is an isomorphism of $S$-modules.
This isomorphism induces an isomorphism of functors
$$
\hat{\gamma}(\alpha,\beta)\, :\;\dual{S\otimes (A,\inv)}_{(\beta\circ\alpha)^{\natural}(I)}\,\xrightarrow{\;\simeq\;}\,
\dual{S\otimes (A,\inv)}_{\beta^{\natural}(U^{-1}I)}\, ,
$$
which is a duality transformation for the identity functor, \ie we have
a duality preserving isomorphism

\smallbreak

$(\id_{\Dfg{S\otimes_{R}A}},\hat{\gamma}(\alpha,\beta))$:
$$
(\Dfg{S\otimes_{R}A},\dual{S\otimes (A,\inv)}_{(\beta\circ\alpha)^{\natural}(I)})
\,\too\, (\Dfg{S\otimes_{R}A},\dual{S\otimes (A,\inv)}_{\beta^{\natural}(U^{-1}I)})\, .
$$
\end{emptythm}

\begin{emptythm}
\label{maintechLem}
{\bf Lemma.}
{\it
Denote $\xi:A\too U^{-1}A$ and $\xi_{1}:U^{-1}A\too S\otimes_{U^{-1}R}U^{-1}A\simeq S\otimes_{R}A$ the natural
$R$-algebra respectively $U^{-1}R$-algebra morphisms. Then
$$
\big(\,(\beta\circ\alpha)_{\ast},\vartheta^{\xi_{1}\circ\xi}\,\big)\, =\,
(\alpha_{\ast},\theta)\circ (\beta_{\ast},\vartheta^{\xi_{1}})\circ (\id_{\Dfg{S\otimes_{R}A}},\hat{\gamma}_{\alpha,\beta})\, .
$$
}

\begin{proof}
The duality preserving functor $\big(\,(\beta\circ\alpha)_{\ast},\vartheta^{\xi_{1}\circ\xi}\,\big)$ on the left
hand side maps $M_{\smb}\in\Dfg{S\otimes_{R}A}$ onto
$$
(\beta\circ\alpha)_{\ast}(M_{\smb})=\alpha_{\ast}(\beta_{\ast}(M_{\smb}))\;\in\,\DcS{A}{\Jideal}\, ,
$$
where the $S\otimes_{R}A$-module $M_{i}$ becomes an $A$-module via the
homomorphism of $R$-algebras $\xi_{1}\circ\xi:A\too S\otimes_{R}A$ for all $i\in\Z$.
The same holds for the duality preserving functor on the right hand side.

\medbreak

Hence we are left to show that the duality transformation of both sides coincide
for all $M_{\smb}\in\Dfg{S\otimes_{R}A}$.

\smallbreak

By the definition of the composition of duality preserving functors,
see~(\ref{dualfunctorCompositionEq}), the duality transformation for the functor
$\alpha_{\ast}\circ\beta_{\ast}\circ\id_{\Dfg{S\otimes A}}$
on the right hand side is given by
$$
\theta_{\beta_{\ast}M}\circ\alpha_{\ast}\big(\,\vartheta^{\xi_{1}}_{M}\circ\beta_{\ast}(\hat{\gamma}(\alpha,\beta)_{M})\,\big)
$$
for $M_{\smb}\in\Dfg{S\otimes_{R}A}$.
We have to show that this is equal to $\vartheta^{\xi_{1}\circ\xi}_{M}$, or equivalently
since $\theta_{\beta_{\ast}M}^{-1}=\eta_{\beta_{\ast}M}$, that
$$
\eta_{\beta_{\ast}M}\circ\vartheta^{\xi_{1}\circ\xi}_{M}\, =\,
\alpha_{\ast}\big(\,\vartheta^{\xi_{1}}_{M}\circ\beta_{\ast}(\hat{\gamma}(\alpha,\beta)_{M})\,\big)\, ,
$$
see~(\ref{ThetaDefEq}).
Now the duality transformation $\vartheta^{\xi_{1}\circ\xi}_{M}$ has in degree~$l$ the $r$-component
$$
\ol{\Hom_{S}(M_{-l-r},\Hom_{R}(S,I_{-r}))}\,\too\,\ol{\Hom_{R}(M_{-l-r},I_{-r})}\, ,
$$
$$
h\,\longmapsto\,\big\{\, m\mapsto\, (-1)^{r}h(m)(1)\,\big\}\, .
$$
Composing this map with the $r$-component in degree~$l$ of $\eta_{\beta_{\ast}M}$ we get
$$
\ol{\Hom_{S}(M_{-l-r},\Hom_{R}(S,I_{-r}))}\,\too\,\ol{\Hom_{U^{-1}R}(M_{-l-r},U^{-1}I_{-r})}
$$
$$
h\,\longmapsto\,\Big\{\, m\mapsto\, (-1)^{r}\iota^{I_{-r}}\big(\,h(m)(1)\,\big)\,\Big\}\, ,
$$
where $\iota^{I_{-r}}:I_{-r}\too U^{-1}I_{-r}$ is the localization morphism.

\smallbreak

But this is the $r$-component in degree~$l$ of
$\alpha_{\ast}\big(\,\vartheta^{\xi_{1}}_{M}\circ\beta_{\ast}(\hat{\gamma}(\alpha,\beta)_{M})\,\big)$.
\end{proof}
\end{emptythm}

\goodbreak
\section{The hermitian Gersten-Witt complex}
\label{GWComplSect}\bigbreak

\begin{emptythm}
\label{FiltsuppSubSect}
{\bf The codimension by support filtration.}
We refer to~\cite{Gi07,Gi09,Gi13} for proofs, details and more information on the hermitian
Gersten-Witt spectral sequence.

\smallbreak

Throughout this section~$X$ denotes a regular and noetherian scheme, and $(\Az,\inv)$
an Azumaya algebra with involution of the first- or second kind over~$X$.

\smallbreak

On~$\DP{\Az}$ we have the filtration by codimension of support:
$$
\DP{\Az}\, =\,\DPF{\Az}{0}\,\supseteq\,\DPF{\Az}{1}\,\supseteq\,\DPF{\Az}{2}\,\supseteq\,\ldots\, ,
$$
where
$$
\DPF{\Az}{p}\, :=\;\Big\{\,\sheaf{F}_{\smb}\,\in\DP{\Az}\,\big|\,\codimension_{X}\supp\sheaf{F}_{\smb}\geq p\,\Big\}
$$
for $p\geq 0$. This is a thick saturated triangulated subcategory of~$\DP{\Az}$ and the duality $\dual{\Az,\inv}$ maps
it into itself. Balmer's~\cite{Ba00} localization sequence gives long exact sequences of Witt groups:
{\small
$$
\xymatrix{
\ldots\too\W^{i}(D^{p}_{\Az})\too\W^{i}(D^{p}_{\Az}/D^{p+1}_{\Az})\xrightarrow{\,\partial\,}
   \W^{i+1}(D^{p+1}_{\Az})\too\W^{i+1}(D^{p}_{\Az})\too\,\ldots\, ,
}
$$
}

\noindent
where we have set $D^{p}_{\Az}\, :=\;\DPF{\Az}{p}$ for all $p\geq 0$,
and the triangular Witt groups are with respect to the duality (induced by) $\dual{\Az,\inv}$.
\end{emptythm}

\begin{emptythm}
\label{GWSpSeqSubSect}
{\bf The hermitian spectral sequence.}
By Massey's method of exact couples we get from the exact sequences above a spectral sequence
$$
E_{1}^{p,q}(\Az,\inv)\, :=\;\W^{p+q}(\DPF{\Az}{p}/\DPF{\Az}{p+1})\, ,
$$
the {\it hermitian Gersten-Witt spectral sequence} of~$(\Az,\inv)$, which converges to the derived hermitian
Witt theory of~$(\Az,\inv)$ if $\dim X<\infty$.

\smallbreak

In~\cite{Gi07,Gi09} it is proven that the odd lines of the hermitian Gersten-Witt spectral sequence are zero,
the lines $E_{1}^{p,4m}(\Az,\inv)$, $m\in\Z$, are all isomorphic to the
{\it hermitian Gersten-Witt complex} of~$(\Az,\inv)$, denoted $\GW_{1}(\Az,\inv)$:
{\footnotesize
$$
\bigoplus\limits_{x\in X^{(0)}}\W_{1}(\Az (x),\inv (x))\too\bigoplus\limits_{x\in X^{(1)}}\W_{1}(\Az (x),\inv (x))\too
\bigoplus\limits_{x\in X^{(2)}}\W_{1}(\Az (x),\inv (x))\too\,\ldots\, ,
$$
}

\smallbreak

\noindent
and the lines $E_{1}^{p,4m+2}(\Az,\inv)$, $m\in\Z$, are all isomorphic to the
{\it skew-hermitian Gersten-Witt complex} of~$(\Az,\inv)$, denoted $\GW_{-1}(\Az,\inv)$:
{\footnotesize
$$
\bigoplus\limits_{x\in X^{(0)}}\W_{-1}(\Az (x),\inv (x))\too\bigoplus\limits_{x\in X^{(1)}}\W_{-1}(\Az (x),\inv (x))\too
\bigoplus\limits_{x\in X^{(2)}}\W_{-1}(\Az (x),\inv (x))\too\,\ldots\, .
$$
}

\smallbreak

\noindent
Here~$X^{(p)}\subseteq X$ denotes the set of points of codimension~$p$ for~$p\geq 0$, and we have set
$(\Az (x),\inv (x)):=k(x)\otimes_{\OXx}(\Az_{x},\inv_{x})$. We consider~$\GW_{\epsilon}(\Az,\inv)$ as a
cohomological complex with $\bigoplus\limits_{x\in X^{(p)}}\W_{\epsilon}(\Az (x),\inv (x))$ in degree~$p$
and denote the $p$th cohomology group of~$\GW_{\epsilon}(\Az,\inv)$ by~$\HM^{p}_{\epsilon}(\Az,\inv)$,
$\epsilon\in\{\pm 1\}$. As usual if~$X=\Spec R$ we use 'affine' notations.
\end{emptythm}

\begin{emptythm}
\label{GWConjSubSect}
{\bf The Gersten conjecture.}
Let now $X=\Spec R$ be an affine scheme associated with a regular integral domain~$R$ of finite Krull dimension
with fraction field~$K$, and set $A:=\Gamma (X,\Az)$.

\smallbreak

The pull-back along the embedding $\iota:R\hookrightarrow K$ induces a homomorphism
$$
\iota^{\ast}\, :\;\W_{\epsilon}(A,\inv)\,\too\, W_{\epsilon}(K\otimes_{R}(A,\inv))
$$
for $\epsilon=\pm 1$. Extending the $\epsilon$-hermitian Gersten-Witt complex by this map on
the left hand side we get the (so called) {\it augmented $\epsilon$-hermitian Gersten-Witt complex}:
{\footnotesize
$$
0\too\W_{\epsilon}(A,\inv)\too\W_{\epsilon}(K\otimes_{R}(A,\inv))\too
\bigoplus\limits_{\hgt q=1}\W_{\epsilon}(k(q)\otimes_{R}(A,\inv))\too\ldots\qquad\qquad\qquad
$$
$$
\qquad\qquad\qquad\qquad\qquad\qquad\qquad\qquad\ldots\too
\bigoplus\limits_{\hgt q=\dim R}\W_{\epsilon}(k(q)\otimes_{R}(A,\inv))\too 0\, ,
$$
}

\smallbreak

\noindent
$\epsilon=\pm 1$. The {\it Gersten conjecture} claims that these complexes are exact if~$R$ is
a regular local ring. By construction this is equivalent to the assertion that the homomorphism
$$
W^{i}(\DPF{A}{p+1},\dual{A,\inv})\,\too\,\W^{i}(\DPF{A}{p},\dual{A,\inv})
$$
is the zero map for all $i\in\Z$ and integers~$p\geq 0$.
\end{emptythm}

\begin{emptythm}
\label{ConsequencesGWConjSubSect}
{\bf A consequence of the Gersten conjecture.}
Given a  regular scheme~$X$ and Azumaya algebra~$(\Az,\inv)$ with involution of the first-
or second kind over~$X$ we denote by $\sheaf{W}^{\epsilon}_{\Az,\inv}$, $\epsilon\in\{\pm 1\}$,
the Zariski sheaf (on~$X$) associated with the presheaf
$$
U\,\longmapsto\,\W_{\epsilon}(\Az|_{U},\inv|_{U})\, ,
$$
where $U\subseteq X$ is an open subscheme. If the Gersten conjecture holds for the Azumaya
algebras with involution of the first- or second kind $(\Az_{x},\inv_{x})$ over~$\OXx$ for all $x\in X$
then the $\epsilon$-hermitian Gersten-Witt complex is a flasque resolution of
$\sheaf{W}^{\epsilon}_{\Az,\inv}$. In particular, we have then
$\HM^{i}_{\Zar}(X,\sheaf{W}^{\epsilon}_{\Az,\inv})\simeq\HM^{i}_{\epsilon}(\Az,\inv)$
for all integers~$i\geq 0$ and all~$\epsilon\in\{\pm 1\}$.

\smallbreak

Another consequence of the Gersten conjecture is the following lemma,
see \eg~\cite[Proof of Cor.\ 7.5]{Gi13} for a proof.
\end{emptythm}

\begin{emptythm}
\label{GWConjLem}
{\bf Lemma.}
{\it
Let~$R$ be a regular local ring, $t\in R$, such that $R/Rt$ is regular as well, and $(A,\inv)$ an Azumaya algebra
with involution of the first- or second kind over~$R$. Assume that the Gersten conjecture holds for $(A,\inv)$.
Then

\smallbreak

\begin{itemize}
\item[(i)]
$\W^{2i+1}(A,\inv)=0$ for all $i\in\Z$; and

\smallbreak

\item[(ii)]
$\HM^{i}_{\Zar}(R_{t},\sheaf{W}_{A,\inv}^{\epsilon}|_{\Spec R_{t}})\, =0$ for all $i\geq 1$ and $\epsilon\in\{\pm 1\}$,
if the Gersten conjecture holds also for $R/Rt\otimes_{R}(A,\inv)$.
\end{itemize}
}
\end{emptythm}

\goodbreak
\section{The main theorem}
\label{GWConjSmoothDVRSect}\bigbreak

\begin{emptythm}
\label{NotationsGWProofSubSect}
Let~$V$ be a field or a discrete valuation ring, $V\too R$ a smooth morphism
of relative dimension~$d$ with~$R$ an integral domain, and~$P$ a prime
ideal of~$R$. Set~$\tilde{R}:=R_{P}$. This is a regular local ring {\it essentially smooth
over~$V$}. Denote by $\gamma:R\too\tilde{R}$ the localization morphism, and by
$\iota:\tilde{R}\hookrightarrow K$ the embedding of~$\tilde{R}$ into its fraction field.

\smallbreak

Let further~$(\tilde{A},\tilde{\inv})$ be an Azumaya algebra with involution of the first- or second kind over~$\tilde{R}$.
Replacing~$R$ by a localization we can assume that $(\tilde{A},\tilde{\inv})=R\otimes_{R}(A,\inv)$ for an Azumaya algebra
with involution of the same kind $(A,\inv)$ over~$R$.

\smallbreak

We denote $I_{\smb}:I_{0}\too I_{-1}\too\ldots\too I_{-\dim R}\in\Dfg{R}$ a minimal injective resolution of~$R$ and
set $\tilde{I}_{\smb}:=I_{\smb\, P}\in\Dfg{\tilde{R}}$.
\end{emptythm}

\begin{emptythm}
\label{mainThm}
{\bf Theorem.}
{\it
The Gersten conjecture holds for the Azumaya algebra with involutions~$(\tilde{A},\tilde{\inv})$ over~$\tilde{R}$.
}

\smallbreak

\noindent
Using the desingularization theorem of Popescu~\cite{Po85,Po86} we get the following more
general case of the Gersten conjecture.

\medbreak

\noindent
{\bf Corollary.}
{\it
Let~$S$ be a regular local ring, which either contains a field, or which is geometrically regular
over a discrete valuation ring. Let further~$(B,\binv)$ be an Azumaya algebra with involution of the
first- or second kind over~$S$. Then the Gersten conjecture holds for~$(B,\binv)$.
}

\begin{proof}
In case~$S$ contains a field it is shown in~\cite[Sect.\ 7]{Gi13} that Theorem~\ref{mainThm}
implies this corollary. Essentially the same arguments work if~$S$ is geometrically
regular over a discrete valuation ring~$V$. We briefly recall the details.

\smallbreak

Let~$f$ be an uniformizer of~$V$. Then since~$S$ is geometrically regular over~$V$ the quotient~$S/Sf$
is regular. It contains the residue field of~$V$, and so the Gersten conjecture holds for $(B',\binv'):=S/Sf\otimes_{S}(B,\binv)$
by the already proven case that the regular local ring contains a field. Analogous since the localization~$S_{f}$
contains the fraction field of~$V$ the Gersten conjecture holds for $(B_{Q},\binv_{Q})$ for all $Q\in\Spec S_{f}$.
Hence, see~\ref{ConsequencesGWConjSubSect}, we have
\begin{equation}
\label{mainCorEq}
\HM_{\Zar}^{i}(S_{f},\sheaf{W}^{\epsilon}_{B_{f},\binv_{f}})\,\simeq\,\HM^{i}_{\epsilon}(B_{f},\binv_{f})
\end{equation}
for all integers~$i\geq 0$ and $\epsilon\in\{\pm 1\}$.

\smallbreak

We use now a consequence of Popescu's desingularization theorem, see~\cite[Cor.\ 1.3]{Sw95}:
The ring~$S$ is a filtered colimit of regular local rings~$S_{\omega}$, $\omega\in\Omega$, which
are essentially smooth over~$V$:
$S=\lim\limits_{{\tiny \begin{array}{c} \too \\ w\in\Omega\end{array}}}S_{\omega}$.
By shrinking the index set~$\Omega$ if necessary we can assume that there exists Azumaya algebras
with involution $(B_{\omega},\binv_{\omega})$ of the same kind as $(B,\binv)$, such that
$(B,\binv)=S\otimes_{S_{\omega}}(B_{\omega},\binv_{\omega})$ for all $\omega\in\Omega$.

\smallbreak

By our main result, Theorem~\ref{mainThm}, the Gersten conjecture holds for $(B_{\omega},\binv_{\omega})$,
$S_{\omega\, Q}\otimes_{S_{\omega}}(B_{\omega},\binv_{\omega})$ for all $Q\in\Spec S_{\omega}$,
and also for $S_{\omega}/S_{\omega}f\otimes_{S_{\omega}}(B_{\omega},\binv_{\omega})$ since
$S_{\omega}/S_{\omega}f$ is essentially smooth over the residue field of~$V$. Therefore by~\ref{ConsequencesGWConjSubSect}
and Lemma~\ref{GWConjLem} we have $\W^{2i+1}(B_{\omega},\binv_{\omega})=0$ for all~$i\in\Z$, and
$$
\HM^{i}_{\epsilon}((B_{\omega})_{f},(\binv_{\omega})_{f})\,\simeq\,
\HM_{\Zar}^{i}((S_{\omega})_{f},\sheaf{W}^{\epsilon}_{B_{\omega},\binv_{\omega}}|_{\Spec S_{\omega\, f}})\; =\, 0
$$
for all $\omega\in\Omega$, $i\geq 1$, and $\epsilon\in\{\pm 1\}$.
Now, see \eg~\cite[Sect.\ 7.1]{Gi13}, we have
$$
\lim\limits_{{\tiny \begin{array}{c} \too \\ w\in\Omega\end{array}}}
\HM_{\Zar}^{i}((S_{\omega})_{f},\sheaf{W}^{\epsilon}_{B_{\omega},\binv_{\omega}}|_{\Spec (S_{\omega})_{f}})\,\simeq\,
\HM_{\Zar}^{i}(S_{f},\sheaf{W}^{\epsilon}_{B,\binv}|_{\Spec S_{f}})\, ,
$$
for all integers~$i\geq 0$ and $\epsilon\in\{\pm 1\}$, and by~\cite[Thm.\ 1.7]{Gi03}
$$
\lim\limits_{{\tiny \begin{array}{c} \too \\ w\in\Omega\end{array}}}\W^{j}(B_{\omega},\binv_{\omega})\,\simeq\,\W^{j}(B,\binv)
$$
for all~$j\in\Z$. We conclude from these considerations taking~(\ref{mainCorEq}) into account that
$$
\W^{2j+1}(B,\binv)\, =\,\HM^{i}_{\epsilon}(B_{f},\binv_{f})\; =\, 0
$$
for all $j\in\Z$, integers~$i\geq 1$, and $\epsilon\in\{\pm 1\}$.

\smallbreak

Since the Gersten conjecture holds for $(B',\binv')=S/Sf\otimes_{S}(B,\binv)$ we also have
$\HM^{i}_{\epsilon}(B',\binv')=0$ for all $i\geq 1$ and $\epsilon\in\{\pm 1\}$

\medbreak

We apply this to the exact cohomology sequence associated with the short exact sequence
of complexes
$$
0\too\GW_{\epsilon}(B',\binv')[-1]\too\GW_{\epsilon}(B,\binv)\too\GW_{\epsilon}(B_{f},\binv_{f})\too 0\, ,
$$
$\epsilon\in\{\pm 1\}$, where $\GW_{\epsilon}(B',\binv',)[-1]$ is the complex $\GW_{\epsilon}(B',\binv')$
shifted by one, \ie starting in degree~$1$, and get $\HM^{i}_{\epsilon}(B,\binv)=0$
for all~$i\geq 2$ and $\epsilon\in\{\pm 1\}$.

\smallbreak

Since~$E_{1}^{p,q}(B,\binv)\Longrightarrow\W^{p+q}(B,\binv)$ and $E_{1}^{p,2m+1}(B,\binv)=0$ for
all $m\in\Z$, see~\ref{GWSpSeqSubSect}, we deduce moreover $\HM^{1}_{\epsilon}(B,\binv)=0$,
and that the natural homomorphism
$$
\W^{1-\epsilon}(B,\binv)\,\too\,\HM^{0}_{\epsilon}(B,\binv)
$$
is an isomorphism for all $\epsilon\in\{\pm\}$. By the main result of Balmer~\cite{Ba01a} we have
$\W_{\epsilon}(B,\binv)\simeq\W^{1-\epsilon}(B,\binv)$, and we are done.
\end{proof}

\medbreak

The proof of Theorem~\ref{mainThm} will follow from the following technical result,
which we prove in Section~\ref{PfGWPfmainLemSect}.
\end{emptythm}

\begin{emptythm}
\label{GWPfmaintechLem}
{\bf Lemma.}
{\it
Let~$t\in R$ be a non zero divisor and non unit, such that $R/Rt$ is flat over~$V$ (which is automatic
if~$V$ is a field). We denote by $\pi:R\too R/Rt=:R'$ the quotient morphism, and set $(A',\inv'):=R'\otimes_{R}(A,\inv)$.

\smallbreak

Let $p\geq 0$ be a natural number and $(M_{\smb},\varphi)$ a $i$-symmetric space in the triangulated
category with duality $(\DfgF{A'}{p},\dual{A',\inv'}_{\pi^{\natural}(I)})$.
Then $\gamma^{\ast}\big(\,\tr_{\pi}(M_{\smb},\varphi)\big)$ is neutral in
$(\DfgF{\tilde{A}}{p},\dual{\tilde{A},\tilde{\inv}}_{\tilde{I}})$.
}

\bigbreak

Before showing that Theorem~\ref{mainThm} follows from this lemma we record a consequence of it.
\end{emptythm}

\begin{emptythm}
\label{GWPfCor}
{\bf Corollary.}
{\it
Let~$t\in\tilde{R}$, such that $\tilde{R}/\tilde{R}t$ is flat over~$V$. Then
$$
\W^{i}(\tilde{A},\tilde{\inv})\,\too\,\W^{i}(\tilde{A}_{t},\tilde{\inv_{t}})
$$
is injective for all $i\in\Z$.
}

\begin{proof}
By Balmer's~\cite{Ba00} localization sequence we have an exact sequence
$$
\W^{i}_{\tilde{R}t}(\tilde{A},\tilde{\inv})\,\too\,\W^{i}(\tilde{A},\tilde{\inv})\,\too\,\W^{i}(\tilde{A}_{t},\tilde{\inv_{t}})\, ,
$$
and so it is enough to show that $\W^{i}_{\tilde{R}t}(\tilde{A},\tilde{\inv})\too\W^{i}(\tilde{A},\tilde{\inv})$ is the
zero map for all $i\in\Z$. By the identification of derived and coherent Witt groups
this is equivalent to show that
$\cws{\tilde{A},\tilde{\inv},\tilde{I}_{\smb}}{\tilde{R}t}{i}\too\cw{\tilde{A},\tilde{\inv},\tilde{I}_{\smb}}{i}$
is trivial for all $i\in\Z$.

\smallbreak

Let~$(\tilde{N}_{\smb},\tilde{\phi})$ be a $i$-symmetric space representing an element of
$\cws{\tilde{A},\tilde{\inv},\tilde{I}_{\smb}}{\tilde{R}t}{i}$. Replacing~$R$ by a localization
we can assume that $t\in R$, $R':=R/Rt$ is flat over~$V$, and
$(\tilde{N}_{\smb},\tilde{\phi})=\gamma^{\ast}(N_{\smb},\phi)$ for some $i$-symmetric
space~$(N_{\smb},\phi)$ in $\DcS{A}{Rt}$ for the duality $\dual{A,\inv}_{I}$, where
$\gamma:R\too\tilde{R}$ is the localization morphism.

\smallbreak

By the d\'evissage theorem~\cite[Thm.\ 5.2]{Gi07} we can assume that
$(N_{\smb},\phi)=\tr_{\pi}(M_{\smb},\varphi)$ for some $i$-symmetric space
in $\Dfg{A'}$ for the duality $\dual{A',\inv'}_{\pi^{\natural}(I)}$, where $\pi:R\too R/Rt$
is the quotient morphism and $(A',\inv'):=R'\otimes_{R}(A,\inv)$.

\smallbreak

Now Lemma~\ref{GWPfmaintechLem} gives
$(\tilde{N}_{\smb},\tilde{\phi})=\gamma^{\ast}(N_{\smb},\phi)=\gamma^{\ast}(\tr_{\pi}(M_{\smb},\varphi))$
is neutral in the triangulated category with duality
$(\DcF{\tilde{A}}{0}\, ,\dual{\tilde{A},\tilde{\inv}}_{\tilde{I}})$. But the $i$th triangular
Witt group of this category with duality is $\cw{\tilde{A},\tilde{\inv},\tilde{I}}{i}$, and we are done.
\end{proof}
\end{emptythm}

\begin{emptythm}
\label{PfmainThmSubSect}
{\bf Proof of Theorem~\ref{mainThm}.}
Modulo some technical details we follow essentially Gillet and Levine~\cite[Proof of Cor.\ 6]{GiLe87}.

\smallbreak

By the construction of the hermitian Gersten-Witt complex we have to show that
$\W^{i}(\DPF{\tilde{A}}{p+1})\too\W^{i}(\DPF{\tilde{A}}{p})$, or equivalently (using the identification
of coherent and derived hermitian Witt groups, see~\ref{der-cohWGrSubSect}), that
$$
\W^{i}(\DfgF{\tilde{A}}{p+1},\dual{\tilde{A},\tilde{\inv}}_{\tilde{I}})\,\too\,
\W^{i}(\DfgF{\tilde{A}}{p},\dual{\tilde{A},\tilde{\inv}}_{\tilde{I}})
$$
is the zero homomorphism for all~$p\geq 0$ and $i\in\Z$.
For this we distinguish the cases $p\geq 1$ and $p=0$ if~$V$ is not a field.

\medbreak

\noindent
{\it Case $p\geq 1$, or $p\geq 0$ and~$V$ is a field.}

\smallbreak

\noindent
Let $\tilde{x}$ be an element of $\W^{i}(\DfgF{\tilde{A}}{p+1})$. Replacing~$\Spec R$ by a
smaller affine neighbourhood of~$P$ if necessary we can assume that $\tilde{x}$ is in the image of
$$
\W^{i}(\DfgF{A}{p+1})\,\too\,\W^{i}(\DfgF{\tilde{A}}{p+1})\, ,
$$
say $\tilde{x}=\gamma^{\ast}(N_{\smb},\phi)$ for some $i$-symmetric space $(N_{\smb},\phi)$ in
$(\DfgF{A}{p+1},\dual{A,\inv}_{I})$. Since the support of~$N_{\smb}$ has codimension~$\geq 2$ if~$V$ is not
a field there exists~$t\in R$ with $R':=R/Rt$ flat over~$V$ and
$\supp N_{\smb}\subseteq\Spec R'$. By d\'evissage, see~\ref{devSubSect}, we have
$(N_{\smb},\phi)=\tr_{\pi}(M_{\smb},\varphi)$ for some $i$-symmetric space
$$
(M_{\smb},\varphi)\quad\mbox{in}\quad
(\DfgF{A'}{p},\dual{A',\inv'}_{\pi^{\natural}(I)})\, ,
$$
where $\pi:R\too R'$ is the quotient morphism, and $(A',\inv')=R'\otimes_{R}(A,\inv)$.

\smallbreak

We conclude by Lemma~\ref{GWPfmaintechLem} above.

\medbreak

\noindent
{\it Case $p=0$ and~$V$ is a discrete valuation ring.}

\smallbreak

\noindent
Let~$f\in V$ be an uniformizer. The ring $\tilde{R}/\tilde{R}f$ is essentially smooth
over the residue field~$V/Vf$ of~$V$ and so a regular local ring. It follows that $P_{0}:=\tilde{R}f$ is
a prime ideal of height one and consequently $\tilde{R}_{P_{0}}$ is a discrete valuation ring.

\smallbreak

By the main result of~\cite{Gi20} the homomorphism
$$
\W^{i}(\tilde{A}_{P_{0}},\tilde{\inv}_{P_{0}})\,\too\,\W^{i}(K\otimes_{\tilde{R}}A,\id_{K}\otimes\,\tilde{\inv})
$$
is injective. On the other hand, $\tilde{R}_{P_{0}}$ is the localization of~$\tilde{R}$ at the multiplicative closed subset
of~$\tilde{R}$ consisting of all $t\in\tilde{R}$ with $\tilde{R}/\tilde{R}t$ flat over~$V$. Hence by Corollary~\ref{GWPfCor}
also $\W^{i}(\tilde{A},\tilde{\inv})\too\W^{i}(\tilde{A}_{P_{0}},\tilde{\inv}_{P_{0}})$ is injective, and therefore
$$
\iota^{\ast}\, :\;\W^{i}(\tilde{A},\tilde{\inv})\,\too\,\W^{i}(K\otimes_{\tilde{R}}\tilde{A},\id_{K}\otimes\,\tilde{\inv})
$$
is a monomorphism for all~$i\in\Z$ as well. In other notations, this means that
$$
\W^{i}(\DPF{\tilde{A}}{0})\,\too\,\W^{i}(\DPF{\tilde{A}}{0}/\DPF{\tilde{A}}{1})\, ,
$$
where the Witt groups are with respect to the duality $\dual{A,\inv}$, respectively with respect to
the by~$\dual{A,\inv}$ induced duality, is injective and therefore by Balmer's~\cite{Ba00}
localization sequence we get that
$$
\W^{i}(\DPF{\tilde{A}}{1})\,\too\,\W^{i}(\DPF{\tilde{A}}{0})
$$
is the zero map for all $i\in\Z$. We are done.
\end{emptythm}

\goodbreak
\section{Proof of Lemma~\ref{GWPfmaintechLem}.}
\label{PfGWPfmainLemSect}\bigbreak

\begin{emptythm}
\label{PreparationPf1}
{\bf Quillen's normalization lemma and a generalization.}
We continue with the notation of the last section, see~\ref{NotationsGWProofSubSect}
as well as Lemma~\ref{GWPfmaintechLem}.

\smallbreak

By Quillen~\cite[\S7, Lem.\ 5.12]{Qu73} if~$V$ is a field, respectively by Gillet-Levine~\cite[Lem.\ 1]{GiLe87}
otherwise, there exists an open immersion $\theta:R\too R_{0}$ with $P\in\Spec R_{0}$
and a smooth morphism $\iota_{0}:\Gamma:=V[T_{1},\ldots ,T_{d-1}]\too R_{0}$ of relative dimension one,
where~$d$ is the relative dimension of~$R$ over~$V$, such that the composition of
this morphism with the quotient map $\pi_{0}: R_{0}\too R_{0}/R_{0}t$ is quasi-finite, respectively finite if~$V$ is a field.

\smallbreak

Using Lemma~\ref{Pushforward-PullbackLem} we can replace~$R$ by~$R_{0}$
to prove Lemma~\ref{GWPfmaintechLem}, and get a commutative (ignoring the morphism
$\tilde{\Delta}$) diagram:
\begin{equation}
\label{Diag1}
\xymatrix{
 & & & & \tilde{R}/\tilde{R}t
\\
\tilde{R}\otimes_{\Gamma}R' \ar[urrrr]^-{\tilde{s}} & &  R\otimes_{\Gamma}R' \ar[ll]^-{\gamma\otimes\id_{R'}}
    & & R' \ar[u]_-{\gamma'} \ar[ll]^-{u}
\\
\\
\tilde{R}\otimes_{\Gamma}R \ar[uu]^-{\id_{\tilde{R}}\otimes\pi} \ar@/^2pc/[dd]^-{\tilde{\Delta}} & & R\otimes_{\Gamma}R \ar[ll]^-{\gamma\otimes\id_{R}}
    \ar[uu]_-{\id_{R}\otimes\pi} & & R \ar[uu]_-{\pi} \ar[ll]^-{p}
\\
\\
\tilde{R} \ar[uu]^-{\tilde{q}} & & R \ar[uu]_-{q} \ar[ll]^-{\gamma} & & \Gamma \ar[uu]_-{\iota} \ar[ll]^-{\iota}  \rlap{\, ,}
}
\end{equation}
where all squares are cartesian, and:

\smallbreak

\begin{itemize}
\item
$\gamma':R'=R/Rt\too\tilde{R}/\tilde{R}t=:\tilde{R}'$ is the localization morphism;

\smallbreak

\item

$\tilde{s}:\tilde{R}\otimes_{\Gamma}R'\too\tilde{R}'$,
$\tilde{r}\otimes x\mapsto\tilde{r}\cdot\gamma' (x)=\tilde{\pi}(\tilde{r})\cdot\gamma' (x)$,
where $\tilde{\pi}:\tilde{R}\too\tilde{R}'=\tilde{R}/\tilde{R}t$ is the quotient morphism;

\smallbreak

\item
$\tilde{\Delta}:\tilde{R}\otimes_{\Gamma}R\too\tilde{R}$, $\tilde{r}\otimes r\mapsto\tilde{r}\cdot\gamma (r)$
is the 'diagonal';

\smallbreak

\item
$p:R\too R\otimes_{\Gamma}R$, $r\mapsto 1\otimes r$, $q:R\too R\otimes_{\Gamma}R$, $r\mapsto r\otimes 1$,
and $\tilde{q}:\tilde{R}\too\tilde{R}\otimes_{\Gamma}R$, $\tilde{r}\mapsto\tilde{r}\otimes 1$
are the 'projections'; and

\smallbreak

\item
$(\id_{\tilde{R}}\otimes\,\pi)\circ\tilde{q}$ is quasi-finite, respectively finite if~$V$ is a field.
\end{itemize}

\medbreak

\noindent
The morphisms~$\iota$, $p$, $q$~and~$\tilde{q}$ are smooth of relative dimension one,
and so~$u$ is also smooth of relative dimension one. Therefore $\tilde{R}\otimes_{\Gamma}R$
is a regular ring of dimension $1+\dim\tilde{R}$. Since $\iota:\Gamma\too R$ is flat the element
$1\otimes t\in\tilde{R}\otimes_{\Gamma}R$ is a non zero divisor, and therefore $\tilde{R}\otimes_{\Gamma}R'$
is a Gorenstein ring of dimension $\dim\tilde{R}$.

\smallbreak

It follows now from~\cite[Chap.\ II, Thm.\ 4.15]{SGA1} that $\tilde{s}:\tilde{R}\otimes_{\Gamma}R'\too\tilde{R}'$
and $\tilde{\Delta}:\tilde{R}\otimes_{\Gamma}R\too\tilde{R}$ are regular embeddings of codimension one.
\end{emptythm}

\begin{emptythm}
\label{PfStrategySubSect}
Set $\tilde{p}:=(\gamma\otimes\,\id_{R})\circ\, p$ and $q_{1}:=(\id_{\tilde{R}}\otimes\,\pi)\circ\tilde{q}$.

\smallbreak

If $\tilde{q}^{\ast}(\tilde{A},\tilde{\inv})\simeq\tilde{p}^{\ast}(A,\inv)$ and~$q_{1}$
is finite we can now finish the proof of Lemma~\ref{GWPfmaintechLem}
as explained in the introduction. However in general~$q_{1}$ is quasi-finite
only and it is possible that $\tilde{q}^{\ast}(\tilde{A},\tilde{\inv})\not\simeq\tilde{p}^{\ast}(A,\inv)$.

\smallbreak

The first obstacle can be resolved using Zariski's main theorem, and for the latter we use a construction
due to Ojanguren and the second named author~\cite[Sects.\ 7 and 8]{OjPa01} to get a smooth
morphism of relative dimension zero $\tilde{R}\otimes_{\Gamma}R\xrightarrow{\kappa}C$,
such that there exists an isomorphism of $C$-algebras with involution
$$
\kappa^{\ast}\big(\tilde{p}^{\ast}(A,\inv))\,\xrightarrow{\;\simeq\;}\,\kappa^{\ast}\big(\tilde{q}^{\ast}(\tilde{A},\tilde{\inv})\big)\, .
$$

\smallbreak

The result of this construction is the following technical lemma.
Before we state the result we note that the algebras with involutions (of the first- or second kind)
$\tilde{q}^{\ast}(\tilde{A},\tilde{\inv})$ and $\tilde{p}^{\ast}(A,\inv)$ become naturally isomorphic
after pull-back along~$\tilde{\Delta}:\tilde{R}\otimes_{\Gamma}R\too\tilde{R}$. More precisely,
we have a natural isomorphism
$$
\rho\, :\;\tilde{\Delta}^{\ast}\big(\,\tilde{p}^{\ast}(A,\inv)\,\big)\,\xrightarrow{\;\simeq\;}\,
\tilde{\Delta}^{\ast}\big(\,\tilde{q}^{\ast}(\tilde{A},\tilde{\inv})\,\big)
$$
of algebras with involutions, which fits into the commutative diagam
\begin{equation}
\label{rhoEq}
\xymatrix{
\tilde{R} \ar[rr] \otimes_{\tilde{R}\otimes_{\Gamma}R}(\tilde{R}\otimes_{\Gamma}R)\otimes_{R}A \ar[rr]^-{\rho}_-{\simeq} \ar[rd]_-{\simeq} & &
     \tilde{R}\otimes_{\tilde{R}\otimes_{\Gamma}R}(\tilde{R}\otimes_{\Gamma}R)\otimes_{\tilde{R}}\tilde{R}\otimes_{R}A \ar[ld]^-{\simeq}
\\
 & \tilde{R}\otimes_{R}A \rlap{\, ,} &
}
\end{equation}
where the diagonal arrows are the natural identifications. Note here that on the left hand side~$R$ acts on $\tilde{R}\otimes_{\Gamma}R$
via the right factor, \ie via $\tilde{p}=(\gamma\otimes\id_{R})\circ p$.
\end{emptythm}

\begin{emptythm}
\label{GeoLem}
{\bf Lemma.}
{\it
There exists a commutative diagram
\begin{equation}
\label{Diag2}
\xymatrix{
 & & & & & & \tilde{R}'
\\
\\
\tilde{C} \ar@/^2pc/[uurrrrrr] ^-{\beta} & & C' \ar@/^1pc/[uurrrr]^-{\beta'} \ar[ll]_-{l} & &
        \tilde{R}\otimes_{\Gamma}R' \ar[ll]_-{\kappa'} \ar[uurr]^-{\tilde{s}}
                 & & R' \ar[uu]_-{\gamma'} \ar[ll]_-{\tilde{u}}
\\
\\
D \ar[uu]^-{\alpha} \ar[uurr]^-{\alpha'} & & C \ar[uu]^-{\pi_{C}} &
      & \tilde{R}\otimes_{\Gamma}R \ar[uu]^-{\id_{\tilde{R}}\otimes\pi} \ar[ll]_-{\kappa} & & R \ar[ll]_-{\tilde{p}} \ar[uu]_-{\pi}
\\
 & & \tilde{R} \ar[u]^-{j} \ar[llu]^-{\delta} \ar[urr]_-{\tilde{q}} \rlap{\, ,} & & & &
}
\end{equation}
where $\tilde{u}=(\gamma\otimes\id_{R'})\circ u$,
$C'=C\otimes_{R}R'=C\otimes_{\tilde{R}\otimes_{\Gamma}R}(\tilde{R}\otimes_{\Gamma}R')$,
and $l:C'\too\tilde{C}:=C'_{Q}$ is the localization homomorphism at
$Q:=(\beta')^{-1}(P\tilde{R}')$ (recall that~$\tilde{R}=R_{P}$ and so $P\tilde{R}$ is the
maximal ideal of~$\tilde{R}$.)

\smallbreak

These rings and morphisms satisfy the following:

\smallbreak

\begin{itemize}
\item[(a)]
$\kappa$ is a smooth morphism of relative dimension zero, and so the same holds for~$\kappa'$;

\smallbreak

\item[(b)]
$C'$ is a Gorenstein ring and $\dim\tilde{C}=\dim\tilde{R}$;

\smallbreak

\item[(c)]
$\beta'$ is a regular immersion of codimension one, and so the kernel of~$\beta$ is generated
by a non unit and non zero divisor;

\smallbreak

\item[(d)]
the by~$\alpha'$ induced morphism of affine schemes $\Spec C'\too\Spec D$
is an open immersion, $\delta$~is a finite morphism, and we have
$$
\tilde{\pi}\, =\,\beta\circ\alpha\circ\delta\, =\, \beta'\circ\pi_{C}\circ j\, =\,
\tilde{s}\circ (\id_{\tilde{R}}\otimes\,\pi)\circ\tilde{q}\, ,
$$
where $\tilde{\pi}:\tilde{R}\too\tilde{R}'=\tilde{R}/\tilde{R}t$ is the quotient morphism;
\smallbreak

\item[(e)]
the morphism~$j=\kappa\circ\tilde{q}$ is smooth of relative dimension one,
and has a splitting $\Delta_{C}:C\too\tilde{R}$ with
$\Delta_{C}\circ\kappa=\tilde{\Delta}$; and

\smallbreak

\item[(f)]
there is an isomorphism of $C$-algebras with involution
$$
\chi\, :\; C\otimes_{R}(A,\inv)\,\xrightarrow{\;\simeq\;}\, C\otimes_{\tilde{R}}(\tilde{A},\tilde{\inv})\, ,
$$
such that $\Delta_{C}^{\ast}(\chi)$ coincides with the natural isomorphism
$$
\Delta_{C}^{\ast}\big(\kappa^{\ast}(\tilde{p}^{\ast}A)\big)\,\xrightarrow{\;\simeq\;}\,\tilde{\Delta}^{\ast}(\tilde{p}^{\ast}A)\,\xrightarrow{\;\rho\;}\,
\tilde{\Delta}^{\ast}(\tilde{q}^{\ast}\tilde{A})\,\xrightarrow{\;\simeq\;}\,\Delta_{C}^{\ast}\big(\kappa^{\ast}(q^{\ast}\tilde{A})\big)\, .
$$
\end{itemize}
}

\begin{proof}
Let $\maxid:=\big(\tilde{s}\circ (\id_{\tilde{R}}\otimes\,\pi)\,\big)^{-1}(P\tilde{R}')$. We have
$\maxid\supseteq\Ker\tilde{\Delta}$, and $\tilde{q}^{-1}(\maxid)=P\tilde{R}$ is the maximal ideal of~$\tilde{R}$.
Hence~$\maxid$ is the unique maximal ideal of $\tilde{R}\otimes_{\Gamma}R$, which contains~$\Ker\tilde{\Delta}$,
and~$\tilde{\Delta}$ factors via the regular local ring $\tilde{S}:=(\tilde{R}\otimes_{\Gamma}R)_{\maxid}$. It follows
\begin{equation}
\label{dimEq0}
\dim\tilde{S}\, =\,\dim (\tilde{R}\otimes_{\Gamma}R)\, =\, 1+\dim\tilde{R}
\end{equation}
since $\maxid$ contains the kernel of
$\tilde{\Delta}:\tilde{R}\otimes_{\Gamma}R\too\tilde{R}$, which is
a regular embedding of codimension one.

\smallbreak

We get a diagram
$$
\xymatrix{
\tilde{S} \ar[rd]_-{\Delta_{\tilde{S}}} & \tilde{R}\otimes_{\Gamma}R \ar[l]_-{\tilde{\iota}} \ar@/^1pc/[d]^-{\tilde{\Delta}}
\\
 & \tilde{R} \ar[u]^-{\tilde{q}} \rlap{\, ,}
}
$$
where $\tilde{\iota}:\tilde{R}\otimes_{\Gamma}R\too\tilde{S}$ is the localization morphism and
$\Delta_{\tilde{S}}\circ\tilde{\iota}=\tilde{\Delta}$. In particular, we have
$\Delta_{\tilde{S}}\circ (\tilde{\iota}\circ\tilde{q})=\id_{\tilde{R}}$, and so there is an isomorphism of
$\tilde{R}$-algebras with involutions
$$
\rho_{\tilde{S}}\, :\;\Delta_{\tilde{S}}^{\ast}\big(\tilde{\iota}^{\ast}(\tilde{p}^{\ast}(A,\inv))\big)\,\xrightarrow{\simeq}\,
\tilde{\Delta}^{\ast}(\tilde{p}^{\ast}(A,\inv))\,\xrightarrow{\;\rho\;}\,\tilde{\Delta}^{\ast}(\tilde{q}^{\ast}(\tilde{A},\tilde{\inv}))
\,\xrightarrow{\simeq}\,\Delta_{\tilde{S}}^{\ast}\big(\tilde{\iota}^{\ast}(\tilde{q}^{\ast}(\tilde{A},\tilde{\inv}))\big)\, .
$$
Now by the theorem~\cite[Prop.\ 7.1]{OjPa01} of Ojanguren and the second named author
there exists a finite \'etale morphism $\tilde{h}:\tilde{S}\too\tilde{\tilde{C}}$, such that

\smallbreak

\begin{itemize}
\item
there is an isomorphism of $\tilde{\tilde{C}}$-algebras with involutions
$$
\tilde{\chi}\, :\; (\tilde{h}\circ\tilde{\iota})^{\ast}\big(\tilde{p}^{\ast}(A,\inv)\big)\,\xrightarrow{\;\simeq\;}\,
(\tilde{h}\circ\tilde{\iota})^{\ast}\big(\tilde{q}^{\ast}(\tilde{A},\tilde{\inv})\big)\, ;
$$
and

\smallbreak

\item
a splitting $\Delta_{\tilde{\tilde{C}}}:\tilde{\tilde{C}}\too\tilde{R}$ of $\tilde{h}\circ\tilde{\iota}\circ\tilde{q}$,
such that $\Delta_{\tilde{\tilde{C}}}^{\ast}(\tilde{\chi})=\rho_{\tilde{S}}$.
\end{itemize}

\smallbreak

\noindent
We can extend these data to an open neighbourhood of the maximal ideal $\maxid$:
There exists $b\in (\tilde{R}\otimes_{\Gamma}R)\setminus\maxid$ and a diagram
$$
\xymatrix{
C \ar[rrrrd]_-{\Delta_{C}} & & (\tilde{R}\otimes_{\Gamma}R)_{b} \ar[ll]_-{h} & & \tilde{R}\otimes_{\Gamma}R \ar[ll]_-{\iota}
    \ar@/_2pc/[llll]_-{\kappa} \ar@/^2pc/[d]^-{\tilde{\Delta}}
\\
 & & & & \tilde{R} \ar[u]^-{\tilde{q}} \rlap{\, ,}
}
$$
where $\iota$ is the localization morphism, $\kappa=h\circ\iota$, $\Delta_{C}\circ\kappa=\tilde{\Delta}$, and such that
there exists an isomorphism of $C$-algebras with involutions
$$
\chi\, :\;\kappa^{\ast}(\tilde{p}^{\ast}(A,\inv))\,\xrightarrow{\;\simeq\:}\,\kappa^{\ast}(\tilde{q}^{\ast}(\tilde{A},\tilde{\inv}))\, ,
$$
such that $\Delta_{C}^{\ast}(\chi)$ coincides with the natural isomorphism of $\tilde{R}$-algebras with involutions
$$
\rho\, :\;\tilde{\Delta}^{\ast}(\tilde{p}^{\ast}(A,\inv))\,\xrightarrow{\;\simeq\;}\,\tilde{\Delta}^{\ast}(\tilde{q}^{\ast}(\tilde{A},\tilde{\inv}))
$$
introduced in~(\ref{rhoEq}). Note that by construction~$\maxid\in\Spec (\tilde{R}\otimes_{\Gamma}R)_{b}$ and therefore
since~$h$ is \'etale and finite we get tacking~(\ref{dimEq0}) into account
\begin{equation}
\label{dimEq}
\dim C\, =\,\dim (\tilde{R}\otimes_{\Gamma}R)_{b}\, =\,\dim (\tilde{R}\otimes_{\Gamma}R)_{\maxid}\, =\,
\dim \tilde{R}\otimes_{\Gamma}R\, =\, 1+\dim\tilde{R}\, .
\end{equation}

\smallbreak

We have constructed the following commutative diagram:
\begin{equation}
\label{CRingEq}
\xymatrix{
 & & \tilde{R}'
\\
C' & & \tilde{R}\otimes_{\Gamma}R' \ar[u]_-{\tilde{s}} \ar[ll]_-{\kappa'}
\\
C \ar[u]^-{\pi_{C}} & & \tilde{R}\otimes_{\Gamma}R \ar[u]_-{\id_{\tilde{R}}\otimes\,\pi} \ar[ll]_-{\kappa}
\\
\tilde{R} \ar[u]^-{j} \ar[rru]_-{\tilde{q}} \rlap{\, ,} &
}
\end{equation}
where $C':=C\otimes_{\tilde{R}\otimes_{\Gamma}R}(\tilde{R}\otimes_{\Gamma}R')$, \ie the
middle square is cartesian, and $j:=\kappa\circ\tilde{q}$, which is smooth of constant relative
dimension one. By construction
$$
\kappa'\, :\;\tilde{R}\otimes_{\Gamma}R'\,\too\, (\tilde{R}\otimes_{\Gamma}R')_{b}\,\too\, C'
$$
is the composition of a localization map followed by a finite \'etale morphism, and so quasi-finite and smooth
of relative dimension zero. It follows that $\pi_{C}\circ j=\kappa'\circ (\id_{\tilde{R}}\otimes\,\pi)\circ\tilde{q}$
is quasi-finite as well since $(\id_{\tilde{R}}\otimes\,\pi)\circ\tilde{q}$ is, and that~$C'$ is a Gorenstein ring,
see~\cite[Cor.\ 3.3.15]{CMR} for the latter claim. Moreover, since the non zero divisor $1\otimes t$ is in~$\maxid$
and $b\not\in\maxid$ we have
$$
\dim (\tilde{R}\otimes_{\Gamma}R')_{b}\, =\,\dim (\tilde{R}\otimes_{\Gamma}R')_{\maxid}\, =\,
\dim (\tilde{R}\otimes_{\Gamma}R)\; -1\, ,
$$
which by~(\ref{dimEq}) implies $\dim (\tilde{R}\otimes_{\Gamma}R')_{b}=\dim\tilde{R}$,
and hence $\dim C'=\dim\tilde{R}$.

\smallbreak

We are done except for the existence of~$\beta'$ and the factorization of the quasi-finite morphism
$\pi_{C}\circ j$. For the later we use a version of Zariski's main theorem, see \eg~\cite[p.\ 42, Cor.\ 2]{ALH}.
By this result the quasi-finite morphism $\pi_{C}\circ j$ factors $\tilde{R}\xrightarrow{\delta}D\xrightarrow{\alpha'} C'$
with~$\delta$ finite and the by~$\alpha'$ induced morphism of affine schemes $\Spec C'\too\Spec D$ an open immersion.

\smallbreak

We are left to show that there exists a regular embedding of codimension one $\beta':C'\too\tilde{R}'$,
such that $\beta'\circ\kappa'=\tilde{s}$. As $C'=C\otimes_{\tilde{R}\otimes_{\Gamma}R}(\tilde{R}\otimes_{\Gamma}R')$
it is for the existence enough to show that there exists a morphism $\beta_{C}:C\too\tilde{R}'$, such that
$$
\beta_{C}\circ\kappa\, =\,\tilde{s}\circ (\id_{\tilde{R}}\otimes\,\pi)\, .
$$

\smallbreak

We claim that $\beta_{C}:=\tilde{\pi}\circ\Delta_{C}$ does the job. In fact,
by construction of~(\ref{CRingEq}) we have
$$
\tilde{\pi}\circ\Delta_{C}\circ\kappa\, =\,\tilde{s}\circ (\id_{\tilde{R}}\otimes\,\pi)\circ\tilde{q}\circ\Delta_{C}\circ\kappa
\, =\,\tilde{s}\circ (\id_{\tilde{R}}\otimes\,\pi)\circ\tilde{q}\circ\tilde{\Delta}\, .
$$
Now observe that for $\tilde{r}\in\tilde{R}$ and~$r\in R$ we have
$$
\begin{array}{r@{\;\; =\;\;}l}
      \tilde{s}\big[\, (\id_{\tilde{R}}\otimes\,\pi)\big(\, (\tilde{q}\circ\tilde{\Delta})(\tilde{r}\otimes r)\,\big)\,\Big] &
              \tilde{s}\big(\, (\id_{\tilde{R}}\otimes\,\pi)( (\tilde{r}\cdot\gamma (r))\otimes 1)\,\big) \\[4mm]
           & \tilde{\pi}(\tilde{r}\cdot\gamma (r)) \\[3mm]
           & \tilde{s}\big(\, (\id_{\tilde{R}}\otimes\,\pi)(\tilde{r}\otimes r)\,\big)\, .
\end{array}
$$
Finally, since~$\kappa'$ is smooth of relative dimension~$0$ and~$\tilde{s}$ is a regular immersion of codimension one,
also~$\beta'$ is a regular immersion of codimension one, see \eg~\cite[Chap.\ IV, Prop.\ 3.9]{RRA}. We are done.
\end{proof}
\end{emptythm}

\begin{emptythm}
\label{PfGWDualComplSubSect}
{\bf The dualizing complexes.}
We have on~$D$ the dualizing complex $E_{\smb}:=\delta^{\natural}(\tilde{I}_{\smb})$,
on~$\tilde{C}$ the dualizing complex $\alpha^{\ast}(E_{\smb})$,  and on~$\tilde{R}'$ the two dualizing complexes
$\tilde{\pi}^{\natural}(\tilde{I}_{\smb})=(\beta\circ\alpha)^{\natural}(E_{\smb})$ and $\beta^{\natural}\big(\alpha^{\ast}(E_{\smb})\big)$,
which are isomorphic to each other, see~\ref{maintechLem2SubSect}.

\smallbreak

Since
$$
\tilde{I}_{\smb}\, :\; \tilde{I}_{0}\too\tilde{I}_{-1}\too\ldots\too \tilde{I}_{-\dim\tilde{R}}\;\;\in\,\Dfg{\tilde{R}}
$$
is a (minimal) injective resolution of~$\tilde{R}$ living in the indicated degrees and
$\dim D=\dim\tilde{R}$ by Lemma~\ref{GeoLem}~(b) and~(d) we know by Lemma~\ref{dualComplLem}
that $\cod{E}(P)=\dim D_{P}$ for all $P\in\Spec D$. Therefore the same holds for the restriction
to the localization $\Spec\tilde{C}$, \ie $\cod{\alpha^{\ast}(E)}(P)=\dim (\tilde{C})_{P}$ for all
$P\in\Spec\tilde{C}$. Since~$\tilde{C}$ is a local Gorenstein ring it follows from the uniqueness of
dualizing complexes, see~\cite[Chap.\ V, Thm.\ 3.1]{RD}, that $\alpha^{\ast}(E_{\smb})$ is
an injective resolution of the $\tilde{C}$-module~$\tilde{C}$.
\end{emptythm}

\begin{emptythm}
\label{2TransferSubSect}
{\bf Two transfer maps.}
Along the morphism $\beta:\tilde{C}\too\tilde{R}'$ we have the following two two duality preserving functors
$$
\big(\,\Dfg{\tilde{A}'},\dual{\tilde{A}',\tilde{\inv'}}_{\beta^{\natural}(\alpha^{\ast}(E))}\,\big)\,\too\,
\big(\,\Dfg{\tilde{C}\otimes_{R'}A'},\dual{\tilde{C}\otimes (A',\inv')}_{\alpha^{\ast}(E)}\,\big)
$$
and
$$
\big(\,\Dfg{\tilde{A}'},\dual{\tilde{A}',\tilde{\inv'}}_{\beta^{\natural}(\alpha^{\ast}(E))}\,\big)\,\too\,
\big(\,\Dfg{\tilde{C}\otimes_{\tilde{R}}\tilde{A}},\dual{\tilde{C}\otimes (\tilde{A},\tilde{\inv})}_{\alpha^{\ast}(E)}\,\big)\, ,
$$
where we have set $(\tilde{A}',\tilde{\inv}'):=\tilde{R}'\otimes_{\tilde{R}}(\tilde{A},\tilde{\inv})$.
These correspond to (\cf~\ref{transferSubSect} for notation)
$$
(\beta,\zeta)\, :\; (\tilde{C},\tilde{C}\otimes_{R'}(A',\inv'))\,\too\, (\tilde{R}',\tilde{R}'\otimes_{R}(A,\inv))\, ,
$$
where~$\zeta$ is the $\tilde{C}$-algebra homomorphism
$$
\tilde{C}\otimes_{R'}A'\, =\,\tilde{C}\otimes_{R'}R'\otimes_{R}A\,\too\,\tilde{R}'\otimes_{\tilde{C}}\tilde{C}\otimes_{R'}R'\otimes_{R}A\,
\xrightarrow{\;\simeq\;}\,\tilde{R}'\otimes_{R}A\, ,
$$
$$
\tilde{c}\otimes r'\otimes a\,\longmapsto\,\big(\beta (\tilde{c})\cdot\gamma' (r')\big)\otimes a\, ,
$$
and
$$
(\beta,\xi)\, :\; (\tilde{C},\tilde{C}\otimes_{\tilde{R}}(\tilde{A},\tilde{\inv}))\,\too\, (\tilde{R}',\tilde{R}'\otimes_{R}(A,\inv))\, ,
$$
where~$\xi$ is the $\tilde{C}$-algebra homomorphism
$$
\tilde{C}\otimes_{\tilde{R}}\tilde{A}\, =\,\tilde{C}\otimes_{\tilde{R}}\tilde{R}\otimes_{R}A\,\too\,
\tilde{R}'\otimes_{\tilde{C}}\tilde{C}\otimes_{\tilde{R}}\tilde{R}\otimes_{R}A\,\xrightarrow{\;\simeq\;}\,
\tilde{R}'\otimes_{R}A\, ,
$$
$$
\tilde{c}\otimes\tilde{r}\otimes a\,\longmapsto\,\big(\beta (\tilde{c})\cdot\tilde{\pi}(\tilde{r})\big)\otimes a\, ,
$$
see~\ref{transferSubSect} for notation.

\smallbreak

By Lemma~\ref{GeoLem}~(f) there exists an isomorphism of $\tilde{C}$-algebras with involutions
$$
\id_{\tilde{C}}\otimes\,\chi\, :\;\tilde{C}\otimes_{C}(C\otimes_{R}A)\,\too\,
\tilde{C}\otimes_{C}(C\otimes_{\tilde{R}}\tilde{R}\otimes_{R}A)\, ,
$$
$$
\tilde{c}\otimes c\otimes a\,\longmapsto\, c\otimes\chi (c\otimes a)\, =\,
\big(\tilde{c}\cdot (l\circ\pi_{C})(c)\big)\otimes\chi (1\otimes a)\, .
$$
Here~$C$ on the left hand side is an $R$-algebra via $R\xrightarrow{\tilde{p}}C$
and on the right hand side it is considered as $\tilde{R}$-algebra via
$\tilde{R}\xrightarrow{j}C$.

\smallbreak

These three morphisms of $\tilde{C}$-algebras with involutions are related as follows.

\medbreak

\noindent
{\bf Lemma.}
{\it
We have a commutative diagram of $\tilde{C}$-algebra morphisms:
$$
\xymatrix{
\tilde{C}\otimes_{C}(C\otimes_{R}A) \ar[rr]^-{\id_{\tilde{C}}\otimes\,\chi} \ar[dd]_-{g} & &
      \tilde{C}\otimes _{C}(C\otimes_{\tilde{R}}\tilde{R}\otimes_{R}A)  \ar[d]^-{\id_{\tilde{C}}\otimes\, g_{1}}
\\
         & & \tilde{C}\otimes_{\tilde{R}}\tilde{R}\otimes_{C}(C\otimes_{\tilde{R}}\tilde{R}\otimes_{R}A)
                  \ar[d]^-{\id_{\tilde{C}}\otimes\, g_{2}}
\\
\tilde{C}\otimes_{R'}R'\otimes_{R}A \ar[rd]_-{\zeta} & &
       \tilde{C}\otimes_{\tilde{R}}(\tilde{R}\otimes_{R}A) \ar[ld]^-{\xi}
\\
 & \tilde{R}'\otimes_{R}A \rlap{\, ,} &
}
$$
where $g,g_{1}$, and~$g_{2}$ are the canonical isomorphisms.
Setting $\ell:=\id_{\tilde{C}}\otimes\, (g_{2}\circ g_{1})$ we have therefore an isometry
$$
\tr_{(\id_{\tilde{C}},g)}\big(\,\tr_{(\beta,\zeta)}(N_{\smb},\psi)\,\big)\,\simeq\,
\tr_{(\id_{\tilde{C}},\,\id_{\tilde{C}}\otimes\chi)}\Big[\,\tr_{(\id_{\tilde{C}},\,\ell)}\big(\,\tr_{(\beta,\xi)}(N_{\smb},\psi)\,\big)\,\Big]
$$
in $(\Dc{\tilde{C}\otimes_{R}A},\dual{\tilde{C}\otimes (A,\inv)}_{\alpha^{\ast}(E)})$
for all $i$-symmetric spaces $(N_{\smb},\psi)$ in the triangulated category with duality
$(\Dc{\tilde{R}'\otimes_{R}A},\dual{\tilde{R}'\otimes (A,\inv)}_{\beta^{\natural}(\alpha^{\ast}(E))}\;)$.
}

\begin{proof}
The last assertion follows from the first, see~\ref{transferSubSect}. To prove that the diagram
commutes we recall first that by Lemma~\ref{GeoLem}~(f) the following diagram commutes
$$
\xymatrix{
\tilde{R}\otimes_{C}(C\otimes_{R}A) \ar[rr]^-{\id_{\tilde{R}}\otimes\,\chi} \ar[rd] & &
    \tilde{R}\otimes_{C}(C\otimes_{\tilde{R}}\tilde{R}\otimes_{R}A) \ar[ld]^-{g_{2}}
\\
 & \tilde{R}\otimes_{R}A \rlap{\, ,} &
}
$$
where the diagonal arrows are the natural isomorphisms. In particular
we have
$$
g_{2}\big(\,(\id_{\tilde{R}}\otimes\,\chi)\big(\tilde{r}\otimes c\otimes a)\big)\, =\,(\tilde{r}\cdot\Delta_{C}(c))\otimes a\, .
$$
Using this we compute for $\tilde{c}\otimes c\otimes a\in\tilde{C}\otimes_{C}C\otimes_{R}A$:

\smallbreak

$\xi\big[\,\ell\big(\, (\id_{\tilde{C}}\otimes\,\chi)(\tilde{c}\otimes (c\otimes a))\,\big)\,\big]$
$$
\begin{array}{r@{\; =\;}l}
    &    \xi\big[\,\ell\big(\, (\tilde{c}\cdot(l\circ\pi_{C})(c))\otimes\,\chi (1_{C}\otimes a)\,\big)\,\big]
\\[4mm]
    & \xi\big[\, (\id_{\tilde{C}}\otimes\, g_{2})\big(\,
          (\tilde{c}\cdot (l\circ\pi_{C})(c))\otimes 1_{\tilde{R}}\otimes\,\chi (1_{C}\otimes a)\,\big)\,\big]
\\[4mm]
    & \xi\Big[\,\big(\tilde{c}\cdot (l\circ\pi_{C})(c)\,\big)\otimes 1_{\tilde{R}}\otimes a\,\Big]
\\[4mm]
    & \beta\Big(\,\tilde{c}\cdot (l\circ\pi_{C})(c)\,\Big)\otimes a\, ,
\end{array}
$$
where we denote for clarity by~$1_{S}$ the one of a ring~$S$.

\smallbreak

On the other hand we have
$$
\zeta\big[\, g\big(\,\tilde{c}\otimes (c\otimes a)\,\big)\,\big]\, =\,
\zeta\Big[\,\big(\tilde{c}\cdot (l\circ\pi_{C})(c)\big)\otimes (1_{R'}\otimes a)\,\Big]\, =\,
\beta\big(\tilde{c}\cdot (l\circ\pi_{C})(c)\big)\otimes a\, ,
$$
hence the lemma.
\end{proof}
\end{emptythm}

\begin{emptythm}
\label{ProofSubSect}
We are now in position to prove Lemma~\ref{GWPfmaintechLem}.
Let for this~$(M_{\smb},\varphi)$ be a $i$-symmetric space in $(\DfgF{A'}{p},\dual{A',\inv'}_{\pi^{\natural}(I)})$.

\smallbreak

The pull-back $\gamma'^{\ast}(M_{\smb},\varphi)$ is a $i$-symmetric space in $\DfgF{\tilde{R}'\otimes_{R}A}{p}$ for the
duality $\dual{\tilde{R}'\otimes_{R}(A,\inv)}_{\tilde{\pi}^{\natural}(\tilde{I})}$, which is isomorphic to the duality
$\dual{\tilde{R}'\otimes_{R}(A,\inv)}_{\beta^{\natural}\big(\alpha^{\ast}(E)\big)}$, see~\ref{PfGWDualComplSubSect}.
For ease of notation we denote by $\gamma'^{\ast}(M_{\smb},\varphi)^{\diamond}$ the space corresponding
to $\gamma'^{\ast}(M_{\smb},\varphi)$ under this isomorphism of triangulated categories, \ie
$$
\gamma'^{\ast}(M_{\smb},\varphi)^{\diamond}\, :=\; \big(\id_{\Dfg{\tilde{R}'\otimes_{R}A}}\, ,\;\hat{\gamma}(\alpha,\beta)\big)_{\ast}
(\gamma'^{\ast}(M_{\smb},\varphi))\, ,
$$
see~\ref{maintechLem2SubSect} for notation.

\smallbreak

We apply now the zero theorem, see~\ref{zeroThmSubSect}. By this result
$$\tr_{(\beta,\zeta)}\big(\,\gamma'^{\ast}(M_{\smb},\varphi)^{\diamond}\,\big)
$$
is a neutral $i$-symmetric space in the triangulated category with duality
$$
\big(\,\DfgF{\tilde{C}\otimes_{R'}A'}{p},\dual{\tilde{C}\otimes_{R'}(A',\inv')}_{\alpha^{\ast}(E)}\,\big)\, .
$$
(Note here that since $\tilde{I}_{\smb}$ is a minimal injective resolution of the $\tilde{R}$-module~$\tilde{R}$ living in
degrees $0,-1,\ldots ,-\dim\tilde{R}$, the complex~$\tilde{\pi}^{\natural}(\tilde{I}_{\smb})$ is a minimal injective
resolution of~$\tilde{R}'$ living in degrees $-1,\ldots , -\dim\tilde{R}$, \cf Example~\ref{dualComplExpl}~(ii).)

\medbreak

By the lemma in~\ref{2TransferSubSect} above this implies

\medbreak

\noindent
{\bf Lemma.}
{\it
The push-forward $\tr_{(\beta,\xi)}\big(\gamma'^{\ast}(M_{\smb},\varphi)^{\diamond}\big)$ is a neutral space in
the triangulated category with duality
$\big(\DfgF{\tilde{C}\otimes_{\tilde{R}}\tilde{A}}{p}\, ,\; \dual{\tilde{C}\otimes_{\tilde{R}}(\tilde{A},\tilde{\inv})}_{\alpha^{\ast}(E)}\big)$.
}

\smallbreak

We compute now:
$$
\begin{array}{r@{\; =\;}l@{\qquad}l}
     \gamma^{\ast}\big(\tr_{\pi}(M_{\smb},\varphi)\big) & \tr_{\tilde{\pi}}\big(\gamma'^{\ast}(M_{\smb},\varphi)\big)
              & \mbox{by Lemma \ref{Pushforward-PullbackLem}} \\[4mm]
         & \tr_{\delta}\big[\,\tr_{\beta\circ\alpha}\big(\gamma'^{\ast}(M_{\smb},\varphi)\big)\,\big] & \\[4mm]
         & \tr_{\delta}\big[\, (\alpha_{\ast},\theta)_{\ast}\big(\,\tr_{(\beta,\xi)}(\gamma'^{\ast}(M_{\smb},\varphi)^{\diamond}\big)\,\big]
              & \mbox{by Lemma~\ref{maintechLem}}
\end{array}
$$
It follows that $\gamma^{\ast}\big(\tr_{\pi}(M_{\smb},\varphi)\big)$
is neutral in $(\DfgF{\tilde{A}}{p},\dual{\tilde{A},\tilde{\inv}}_{\tilde{I}})$
since the space $\tr_{(\beta,\xi)}\big(\,\gamma'^{\ast}(M,\varphi)^{\diamond}\,\big)$ is neutral in
$\big(\DfgF{\tilde{C}\otimes_{\tilde{R}}\tilde{A}}{p}\, ,\; \dual{\tilde{C}\otimes_{\tilde{R}}(\tilde{A},\tilde{\inv})}_{\alpha^{\ast}(E)}\big)$
by the lemma above. We are done.
\end{emptythm}

\begin{emptythm}
\label{wrongGWPfRem}
{\bf Remark.}
In the article~\cite{Gi13} by the first named author it was not observed that
if $(\tilde{A},\tilde{\inv})$ is not extended from the base ring then the $R\otimes_{\Gamma}R$-algebras
with involutions $p^{\ast}(A,\inv)$ and~$q^{\ast}(\tilde{A},\tilde{\inv})$ are not necessarily isomorphic (using the notation
of~(\ref{Diag1})). Hence~\cite{Gi13} proves the Gersten conjecture only in the constant case,
\ie in case~$\tilde{R}$ is a regular local ring which contains a field~$V$ and $(\tilde{A},\tilde{\inv})$
is extended from~$V$.
\end{emptythm}

\bibliographystyle{amsalpha}

\begin{thebibliography}{CTHK96}

\smallbreak

\bibitem{Ba00}
P.\ Balmer, \textsl{Triangular Witt groups I: The 12-term localization exact sequence}, $K$-Theory
{\bfseries 19} (2000), 311--363.

\smallbreak

\bibitem{Ba01a}
P.\ Balmer, \textsl{Triangular Witt groups II: From usual to derived}, Math.\ Z.\ {\bfseries 236}
(2001), 351--382.

%

\smallbreak

\bibitem{BaGiPaWa02}
P.\ Balmer, S.\ Gille, I.\ Panin, C.\ Walter, \textsl{The Gersten conjecture for Witt groups in the equicharacteristic
case}, Doc.\ Math.\ {\bfseries 7} (2002), 203--217.

%

\smallbreak

\bibitem{BaWa02}
P.\ Balmer, C.\ Walter, \textsl{A Gersten-Witt spectral sequence for regular schemes}, Ann.\ Sci.\
\'Ecole Norm.\ Sup.\ (4) {\bfseries 35} (2002), 127--152.

\smallbreak

\bibitem{BaFiPa19}
E.\ Bayer-Fluckiger, U.\ First, R.\ Parimala, \textsl{On the Grothendieck-Serre conjecture for classical groups},
Preprint 2019, {\ttfamily arXiv:1911:02518}.

%

\smallbreak

\bibitem{CMR}
W.\ Bruns, J.\ Herzog, \textsl{Cohen-Macaulay rings}, Cambridge Studies in Advanced Mathematics,
{\bfseries 39}, Cambridge University Press, Cambridge, 1993.

%

%

%

\smallbreak

\bibitem{RRA}
W.\ Fulton, S.\ Lang, \textsl{Riemann-Roch algebra}, Grundlehren der mathematischen Wissenschaften, {\bfseries 277},
Springer-Verlag, New York, 1985.

%

\smallbreak

\bibitem{Gi02}
S.\ Gille, \textsl{On Witt groups with support}, Math.\ Ann.\ {\bfseries 322} (2002), 103--137.

\smallbreak

\bibitem{Gi03}
S.\ Gille, \textsl{Homotopy invariance of coherent Witt groups}, Math.\ Z.\ {\bfseries 244} (2003),
211--233.

\smallbreak

\bibitem{Gi07}
S.\ Gille, \textsl{A Gersten-Witt complex for hermitian Witt groups of coherent algebras over schemes I:
Involution of the first kind}, Compos.\ Math.\ {\bfseries 143} (2007), 271--289.

\smallbreak

\bibitem{Gi09}
S.\ Gille, \textsl{A Gersten-Witt complex for hermitian Witt groups of coherent algebras over schemes II:
Involution of the second kind}, J.\ $K$-Theory {\bfseries 4} (2009), 347--377.

\smallbreak

\bibitem{Gi13}
S.\ Gille, \textsl{On coherent hermitian Witt groups},  Manuscripta Math.\ {\bfseries 141} (2013), 423--446.

\smallbreak

\bibitem{Gi20}
S.\ Gille, \textsl{A hermitian analog of a quadratic form theorem of Springer},  Manuscripta Math.\ {\bfseries 163}
(2020), 125--163.

\smallbreak

\bibitem{GiLe87}
H.\ Gillet, M.\ Levine, \textsl{The relative form of Gersten's conjecture over a discrete valuation ring: the smooth case},
J.\ Pure Appl.\ Algebra {\bfseries 46} (1987), 59--71.

\smallbreak

\bibitem{SGA1}
A.\ Grothendieck, \textsl{Rev\^etements \'etales et groupe fondamental}, S\'eminaire de G\'eom\'etrie
Alg\'ebrique du Bois Marie 1960Ð1961 (SGA 1). Dirig\'e par Alexandre Grothendieck. Augment\'e de deux expos\'es de
M.\ Raynaud, Lecture Notes in Math. {\bfseries 224}, Springer-Verlag, Berlin-New York, 1971.

\smallbreak

\bibitem{RD}
R.\ Hartshorne, \textsl{Residues and duality}, Lecture Notes in Math.\ {\bfseries 20}, Springer-Verlag, Berlin-New York 1966.

%

%

%

%

%

\smallbreak

\bibitem{OjPa01}
M.\ Ojanguren, I.\ Panin, \textsl{Rationally trivial Hermitian spaces are locally trivial}, Math.\ Z.\
{\bfseries 237} (2001), 181--198.

%

\smallbreak

\bibitem{Po85}
D.\ Popescu, \textsl{General N\'eron desingularization}, Nagoya Math.\ J.\ {\bfseries 100} (1985), 97--126.

\smallbreak

\bibitem{Po86}
D.\ Popescu, \textsl{General N\'eron desingularization and approximation}, Nagoya Math.\ J.\ {\bfseries 104} (1986), 85--115.

\smallbreak

\bibitem{Qu73}
D.\ Quillen, \textsl{Higher algebraic K-theory. I}, in Algebraic K-theory, I: Higher K-theories (Proc.\ Conf., Battelle Memorial
Inst., Seattle, Wash., 1972), pp. 85--147. Lecture Notes in Math.\ {\bfseries 341}, Springer-Verlag, Berlin 1973.

\smallbreak

\bibitem{ALH}
M.\ Raynaud, \textsl{Anneaux locaux hens\'eliens}, Lecture Notes in Math.\ {\bfseries 169}, Springer-Verlag, Berlin-New York 1970.

\smallbreak

\bibitem{Sw95}
R.\ Swan, \textsl{N\'eron-Popescu desingularization}, Algebra and geometry (Taipei, 1995), 135--192,
Lect.\ Algebra Geom., {\bfseries 2}, Int.\ Press, Cambridge, MA, 1998.

\smallbreak

\bibitem{Ve96}
J.\ Verdier, \textsl{Des cat\'egories d\'eriv\'ees des cat\'egories ab\'eliennes}, Ast\'eriques {\bfseries 239} (1996).



\end{thebibliography}

\end{document}